\documentstyle[11pt,leqno]{article}

\newfont{\Bb}{msbm10 scaled\magstephalf}

\begin{document}
\noindent

\begin{center}
  {\LARGE  Gromov-Witten invariants and rigidity of Hamiltonian loops with
 compact support  on  noncompact symplectic manifolds}
\end{center}

  \noindent
  \begin{center}
    {\large  Guangcun Lu}\footnote{Supported by the NNSF
   19971045 and ETPME of China.}\\[5pt]
       Department of Mathematics, Beijing Normal University\\
               Beijing 100875, P. R. China\\[5pt]  \vspace{-2mm}
               (E-mail: gclu@bnu.edu.cn)\\[3pt]

              Preliminary version\hspace{2mm}October 1, 1997\\[5pt]
              \vspace{-4mm}
              Revised version\hspace{4mm} May 1, 1999  \\[5pt]
              Final revised version  November 13, 2001\\[5pt]
              \end{center}

\begin{abstract} In this paper the Gromov-Witten invariants on a class of
noncompact symplectic manifolds are defined by combining Ruan-Tian's method
with that of McDuff-Salamon.  The main point of the arguments is to introduce
a method dealing with the transversality problems in the case of noncompact
manifolds. Moreover, the techniques are also used to study the topological
rigidity of Hamiltonian loops with compact support on a class of noncompact
symplectic manifolds.
 \end{abstract}



          \section{Introduction}

Since Gromov introduced his celebrated pseudo-holomorphic theory
on symplectic manifolds in 80's([Gr]), many important questions in
symplectic geometry and related fields have been solved. In
particular, Witten [W1, W2] pointed out that  Gromov's study of
the moduli space of holomorphic curves could be used in principle
to describe correlation functions in the topological quantum field
theory. The moduli spaces of holomorphic spheres were used by Ruan
to define certain symplectic invariants of semi-positive
symplectic manifolds([R1]). In the semi-positive closed symplectic
manifolds the more general Gromov-Witten invariants of any genus,
including so called  mixed invariants, were constructed in [RT1]
and later [RT2] and thus they gave the first rigorous mathematics
theory of quantum cohomology. This forms a solid mathematical
basis for the topological sigma model. In addition, they also
applied these invariants to the Mirror Symmetry Conjecture, the
Enumerative Geometry and Symplectic topology. It should be noted
that in this case their mixed invariants are of integral values.

On the other hand the Gromov-Witten invariants were studied axiomatically in [KM].
More recently, the Gromov-Witten invariants  for any  projective  manifolds(cf.[LT1])
and any closed symplectic manifolds were defined(cf.[FO][LT2][R3][Sie]).

It has been expected that the Gromov-Witten invariants should also be defined
for noncompact symplectic manifolds and  families of symplectic manifolds(cf.[K1, pp.364]).
 In fact, the latter was carried out in [L][R3]. Roughly saying,
if $p:Y\to M$ is an oriented fiber bundle such that the fiber $X$ and
the base $M$ are smooth, compact, oriented manifolds( which implies that
$Y$ is also such a manifold), and $\omega$ is a closed $2$-form on $Y$
such that $\omega$ restricts to a symplectic form over each fiber, then
 $Y$ can be viewed as a family of symplectic manifolds and the
 Gromov-Witten invariants over $Y$ are defined in [R3].
However, for noncompact symplectic manifolds $(V,\omega)$ how the
Gromov-invariants over them  should be defined, we so far do not
see it in the literatures. Generally speaking, the key points in
many applications of the Gromov's pseudo-holomorphic curve theory
are  the compactness problems. On the closed symplectic manifolds
one have obtained very good results(cf. [Gr], [RT1], [PW], [Ye]).
For the general noncompact symplectic manifolds( even  without
boundary) these problems become very complicated. In this paper we
define the Gromov-Witten invariants on a class of special
noncompact symplectic manifolds------semi-positive geometrically
bounded one. Precisely speaking, we generalize the main results in
[RT1] to this class of symplectic manifolds. The notion of
geometrically bounded (abb. g. bounded) symplectic manifolds was
first appeared in [Gr]. This kind of manifolds has many nice
properties so that many results on closed symplectic manifolds can
be extended on them in some reasonable ways( see \S2).

However, since  $V$ is  noncompact, for every integer $m\ge 1$ the
Banach manifolds ${\cal J}^m_\tau$ consisting of all $C^m$-smooth
$\omega$-tame almost complex structures on $(V,\omega)$ and the
group ${\rm Diff}^m(V)$ of all $C^m$-diffeomorphisms on $V$ are
not separable, and thus neither are some correspondent moduli
spaces separable. Hence it is difficult using Sard-Smale theorem
in many transversality arguments. One may wish to use its
generalization version due to Quinn  to replace it. But this
requires the Fredholm map considered to be proper or
$\sigma$-proper. Under our case it can not be satisfied.  On the
other hand, for a given $J\in {\cal J}^m_\tau(M,\omega)$ the space
$C^m({\rm T}_J)$ of all $C^m$-sections does not gives rise to a
local model for the space  ${\cal J}^m_\tau(V,\omega)$ via
$Y\mapsto J{\rm exp}(-JY)$. To see  this point, note that
 $J\in{\cal J}^m_\tau(V,\omega)$ only means
 $\omega(\xi,J(p)\xi)>0$ for every $p\in V$ and $\xi\in T_pV\setminus\{0\}$
 and from $\|Y\|_{C^m}<\delta$ it does not follow that $\|Y\|_{C^0}<\eta$
 which is an arbitrary given positive number smaller than $\delta$.
 Thus even if for every  $p\in V$ and $\xi\in T_pV\setminus\{0\}$ we can
 obtain  $\omega(\xi,J(p){\rm exp}(-J(p)Y(p))\xi)>0$ as $|Y(p)|$
 sufficiently small,  but due to the noncompactness of $V$ one can
 not derive that for a given  smooth nowhere null vector field $\zeta$
 on $V$, $\omega(\zeta(p),J(p){\rm exp}(-J(p)Y(p)) \zeta(p))$ is
 more than zero at all points $p\in V$ whether $\|Y\|_{C^m}$ is small.
In order to overcome these difficulties we construct suitable
separable Banach manifolds to replace the Banach manifolds chosen
naturally in the case of compact manifolds. In \S2 and \S4 these
techniques are all  used. The method may probably applied to
generalize other results on compact manifolds in symplectic
geometry and Seiberg-Witten invariants theory to noncompact
manifolds.

In our case replacing $H^*(V,\mbox{\Bb Z})$ by
 $H_*(V,\mbox{\Bb Z})$ the homology  we can
 show that there is an quantum ring structure on it.
 In contrast to the case of closed symplectic manifolds
 it seem to be very hard to use the recent techniques developed by
 [FO][LT][R3][Sie] to define the Gromov-Witten invariants on all noncompact compact
g.bounded symplectic manifolds  because of the
 technical  difficulties.

 Inspired by Seidel's work [Se1] the quantum homology is also used to study topological
 rigidity of Hamiltonian loops by F.Lalonde, D. McDuff and L. Polterovich in [LMP].
 Precisely speaking, they proved that if $\omega_1$ and $\omega_2$ are two symplectic
 forms satisfying certain monotonicity assumptions on a closed manifold $M$ then
 every loop $\phi=\{\phi_t\}_{0\le t\le 1}$ in the group
 ${\rm Ham}(M,\omega_1)\cap{\rm Symp}(M,\omega_2)$ can be homotoped in
 ${\rm Symp}(M,\omega_2)$ to a loop in ${\rm Ham}(M,\omega_2)$.
 Combing their ideas with our techniques together we generalize their results
 to the case of the Hamiltonian loops with compact support on a class of
 noncompact g.bounded symplectic manifolds in Corollary 6.2. Moreover,
 as a consequence the corresponding result on compact symplectic manifolds
 with contact type boundary is also obtained in Corollary 6.3. The main points
 of the arguments are to construct a kind of suitable closed two-forms on the
 Hamiltonian fibre bundle over $S^2$ with noncompact g.bounded symplectic manifolds
 as a fibre to replace the unique coupling class whose top power vanishes
 so that the composition rule may be obtained.

   The arrangements of  this paper are as follows.
  In \S2 we give some basic
 definitions  and   lemmas in geometrically bounded symplectic manifolds,
 and specially a new technique on transversality arguments. In \S3 we
 generalized the results of transversality and compactness  to our case.
 Since the arguments are similar we only give the necessary  improvements.
 The Gromov-Witten invariants are defined in \S4.   As a consequence
 we also define the Gromov-Witten invariants of compact symplectic
 manifolds with contact type boundary in \S5. In \S6 the study of
 the topological rigidity of Hamiltonian loops with compact support
 on noncompact g.bounded symplectic manifolds with the weaker
 semi-positivity assumptions is given. In final Appendix a theorem
 which characterizes the Hamiltonian symplectomorphisms on a compact
 symplectic manifold with contact type boundary in terms of the
 flux homomorphism is provided.

{\bf Acknowledgements.}\hspace{2mm} This revised version including
the topological rigidity of Hamiltonian loops was accomplished
during the author's visit at IHES. The author would like to
express his thanks to Professor J.Bourguignon for his invitation
and hospitality. He thanks very much Dr. Paul Seidel for carefully
checking the original version, correcting many mistakes and
stimulating conversations. He is also grateful to Professors Gang
Tian and Yongbin Ruan for their kind explanations to me with
respect to their papers. Finally, I wish to thank the referee for
telling me another possible simpler method to achieve the
transversality in the definition of our Gromov-Witten invariants.

\section{Definitions and Lemmas }

In this section we give some necessary technical lemmas. Notice first
that the following conclusions in  Riemannian geometry  are some easy
exercises.

\noindent{\bf Lemma 2.1.}\hspace{2mm}{\it Let $(M, g)$ be a
Riemannian manifold with injectivity radius $i(M, g)>0$. Then it
is complete and for any compact subsets $K$ in $M$ and arbitrary
$\varepsilon>0$,
$$K_{\varepsilon}=\{p\in M: d_g(p, K)\le\varepsilon\}$$
is compact. Here $d_g$ denote the distance induced by $g$.}\vspace{2mm}

\noindent{\bf Lemma 2.2.}\hspace{2mm}{\it For the product
Riemannian manifold $(M, g)=(M_1, g_1)\times(M_2, g_2)$ we have
\begin{description}
\item[(i)]$i(M, g)=\min\{i(M_1, g_1), i(M_2, g_2)\}$;
\item[(ii)]$\forall (m_1, m_2)\in M$, $u=(u_1, u_2)$, and
$v=(v_1, v_2)\in T_mM$ it holds that
$$K_g(\Pi_m)=\frac{1}{4}\biggl(K_{g_1}(\Pi_{m_1})+K_{g_2}(\Pi_{m_2})\biggr),$$
where $\Pi_m=span\{u,v\},\, \Pi_{m_1}=span\{u_1, v_1\}$ and\,
$\Pi_{m_2}=span\{u_2, v_2\}$.
\end{description}}\vspace{2mm}

Next let us recall the following definition(cf. [ALP] [Gr] [Sik]).

\noindent{\bf Definition 2.3.}\hspace{2mm}{\it Let $(V, \omega)$
be a symplectic manifold without boundary. Call it {\rm
geometrically bounded} if there exists  an almost complex
structure $J$ and a complete Riemannian metric $g$ on $V$ such
that the following properties are satisfied:
\begin{description}
\item[$1^{\circ}$] $J$ is uniformly tamed by $\omega$, that is, there exist
strictly positive constants $\alpha_0$ and $\beta_0$ such that
$$\omega(X, JX)\ge\alpha_0\|X\|^2_g\quad {\rm and}\quad
|\omega(X, Y)|\le\beta_0\|X\|_g\|Y\|_g$$
for all $X, Y\in TV$;
\item[$2^{\circ}$]  The sectional curvature $K_g\le C_0$(a positive constant)
and the injectivity  radius  $i(V, g)>0$.
\end{description}}

\noindent{\bf Remark 2.4.}\hspace{2mm} By Lemma 2.1 we know that
the requirement of the completeness for $g$ in Definition 2.3 is
not necessary since this is actually contained in the condition
$2^{\circ}$.

Clearly the closed symplectic manifolds are s.g. bounded, a
product of two g. bounded symplectic manifolds is also such
manifold. One can easily prove that every symplectic covering
manifold of a g. bounded symplectic manifold and every symplectic
manifold without boundary which is isomorphic at infinity  to the
symplectization of a closed contact manifold are g. bounded. In
[Lu2] we have proved that the cotangent bundles with respect to
any twisted symplectic structures on it are g.bounded. In
addition, one also should notice that any geometrically bounded
symplectic manifolds are the tame almost complex manifolds in the
sense of J.C.Sikorav (see[Sik]).

Given a closed Riemann surface $\Sigma$ with the complex structure $j$ and
$J\in{\cal J}_\tau(V,\omega)$ we denote by $\overline{Hom}_J(T\Sigma,TV)$
the space of the smooth sections of the bundle of anti-J-linear homomorphisms
from $T\Sigma$ to $TV$ over $\Sigma\times V$. Its element $\nu$ is called
the inhomogeneous term. Recall that a smooth map $f:\Sigma\to V$ is called
$(J,\nu)$-map if for any $z\in\Sigma$,
$$\bar\partial_Jf(z)=df(z)+J(f(z))\circ df(z)\circ j(z)=\nu(z,f(z)).$$
In the following we only consider the inhomogeneous term $\nu$ satisfying
\begin{equation}
{\rm Sup}_{(z,p)\in\Sigma\times V}\|\nu(z,p)\|_{{\cal L}(T_z\Sigma,T_pV)}
<+\infty
\end{equation}
where the norm in ${\cal L}(T_z\Sigma,T_pV)$ is with respect to
$g$ and the Riemannian metric on $\Sigma$ induced from $j$ and
some area form. Notice that any two area forms on $\Sigma$ are
proportional. The above condition is independent of the concrete
choice of the compatible area forms.

\noindent{\bf Lemma 2.5.}\hspace{2mm}{\it Let $(V,\omega,g,J)$ be
as above Definition 2.3, and $\sigma$ an area form on $\Sigma$
compatible with $j$,$\tau=\sigma\circ(id\times j)$. Then for $N>0$
sufficiently large $(\Sigma\times V, \tilde\omega,\tau\oplus
g,\tilde J)$ is also a g.bounded symplectic manifold. Here
$\tilde\omega=N\tau\times\omega$ and $\tilde
J(z,p):T_{(z,p)}(\Sigma\times V)\to T_{(z,p)}(\Sigma\times V)$ is
given by
$$(X_1, X_2)\mapsto (j(z)X_1, J(p)\nu(z,p)(X_1)+J(p)X_2).$$}

The proof of this lemma is an easy exercise. In fact, one can choose
$(\alpha_1,\beta_1)$ to replace $(\alpha_0,\beta_0)$. Here $\alpha_1=
\alpha_0/2$, $\beta_1=2\beta_0+\alpha_0 +\Gamma^4\beta_0^4/\alpha_0\eta$
and
$${\rm Sup}_{(z,p)\in\Sigma\times V}\|\nu(z,p)\|_{{\cal L}(T_z\Sigma,T_pV)}
\le\Gamma<+\infty,\qquad N\ge\frac{\alpha_0}{2}+\frac{\Gamma^4\beta_0^4}{
2\alpha_0^2}.$$

\noindent{\bf Proposition 2.6.}\hspace{2mm}{\it Under assumptions
of Lemma 2.5, if $K\subset V$ is a compact subset and $u:\Sigma\to
V$ a smooth $(J,\nu)$-map representing $A\in H_2(V,\mbox{\Bb Z})$
and intersecting with $K$, then
$${\rm Im}(f)\subset K_{\rho_0},$$
where $\rho_0=\rho_0(\alpha_0,\beta_0,C_0, i(V,g),j, J,\nu, A,\sigma).$}

\noindent{\it Proof}.\hspace{2mm}Write $W=\Sigma\times V$ and
$\bar u:\Sigma\to W, z\mapsto (z,u(z))$. Then $\bar u$ is $\tilde
J$- holomorphic and its image can intersect with
$\widehat{K}:=\Sigma\times K$
 if and only if
the image of $u$ is intersecting with $K$. Combing this with the taming
property we can estimate  its area with respect to the metric
$\tau\oplus g$  as follows:
\begin{eqnarray*}
{\rm Area}_{\tau\oplus g}(\bar u(\Sigma))
&\le&\frac{1}{\alpha_1}\int_{\Sigma}\bar u^*\tilde\omega\\
 &=&\frac{1}{\alpha_1}\int_{\Sigma}u^*\omega +\frac{N}{\alpha_1}
 \int_{\Sigma}\sigma\\
 &=&\frac{1}{\alpha_1}<\omega,A> +\frac{N}{\alpha_1}\int_{\Sigma}\sigma.
 \end{eqnarray*}
Now, by Lemma 2.5 we have
\begin{eqnarray}
\tilde \omega((X_1, X_2), \tilde J_{\lambda}(X_1, X_2))&\ge &
\frac{\alpha_0}{2}\|(X_1,X_2)\|^2_{\tau_0\oplus g}\\
|\tilde\omega((X_1,X_2),(Y_1,Y_2))|&\le &
\beta_1\|(X_1,X_2)\|_{\tau\oplus g}\|(Y_1,Y_2)\|_{\tau\oplus g}
\end{eqnarray}
for every $(z,p)\in W$ and $X=(X_1, X_2),Y=(Y_1,Y_2)\in T_{(z,p)}W$.
Here $\|(X_1,X_2)\|^2_{\tau\oplus g}=\|X_1\|^2_{\tau}+\|X_2\|^2_g$.
Moreover, by lemma 2.2 the sectional curvature and injectivity radius of
$(\Sigma\times V, \tau\oplus g)$ satisfy
\begin{equation}
K_{\tau\oplus g}\le\frac{1}{4}(1+C_0)
\end{equation}
and
\begin{equation}
i(\Sigma\times V,\tau\oplus g)={\rm min}\{i(\Sigma,\tau),i(V,g)\},
\end{equation}
respectively.
Next, according to the comments below Definition 4.1.1 in [Sik],
in our case  we may take $C_1=1/\pi,\; C_2=\beta_1/\alpha_1,
\;\omega_x\equiv\tilde\omega/\beta_1$ and
$r_0=\min(i(W,\tau\oplus g), 2\pi/\sqrt{1+C_0})$ such that the
following monotonicity holds:

{\it For  a compact Riemannian surface $S$ with boundary and
$\tilde J$-holomorphic map $f: S\to W$, if $f(S)\subset  B(x, r)\subset W$,
$f(\partial  S)\subset\partial B$ and $x\in f(S)$  for some $r\le r_0$,
then
\begin{equation}
{\rm Area}_g(f(S))\ge\frac{\pi\alpha_1}{4\beta_1}r^2.
\end{equation}}
From these and the proof of Proposition 4.41 in [Sik] it follows
that
$${\rm Im}(\bar u)\subset U(\widehat{K},C_6{\rm Area}(Im(\bar u))),$$
where $\widehat{K}=\Sigma\times K$ and $C_6=4C_1C_2/r_0=
4\beta_1/\pi\alpha_1r_0$. Using the argument below Lemma 2.5 and an
easy computation we can get
$$C_6=\frac{4\alpha_0^3\beta_0+2\alpha_0^4+2\beta_0^4\Gamma^4}
{\pi\alpha_0^4\min(i(\Sigma,\tau),i(V,g),2\pi/\sqrt{1+C_0})}.$$
Notice that $${\rm Area}(\bar u(\Sigma))\le\frac{2}{\alpha_0}<\omega,A>+
\frac{2N}{\alpha_0}\int_{\Sigma}\sigma,$$
and we can choose
\begin{equation}
\Gamma={\rm Sup}_{(z,p)\in\Sigma\times V}\|\nu(z,p)\|_{(\tau,g)},\quad
{\rm and}\quad N=\frac{\alpha_0}{2}+\frac{\Gamma^4\beta_0^4}{2\alpha_0^3}.
\end{equation}
Therefore we can find a positive number
\begin{equation}
\rho=\rho(\alpha_0,\beta_0,C_0, i(V,g),j, J,\nu, A,\sigma)
\end{equation}
such that
\begin{equation} {\rm Im}(\bar u)\subset\widehat{K}_\rho.
\end{equation}
Projecting on $V$ we can complete the proof of Proposition 2.6.
\hfill$\Box$\vspace{2mm}

As pointed out in Introduction, generally speaking, on the noncompact manifold
$V$ for a given $J\in{\cal J}^m_\tau(V,\omega)$ and an arbitrary small
 positive number $\delta>0$ there may exist a $C^m$-smooth section $Y$ of the bundle
${\rm T}_J\to V$ such that $\|Y\|_{C^m}<\delta$, but  $J{\rm exp}(-JY)\notin
 {\cal J}^m_\tau(V,\omega)$.
 But for some noncompact symplectic manifolds we can prove:

\noindent{\bf Lemma 2.7.}\hspace{2mm}{\it For a given $J_0\in{\cal
J}^m_\tau(V,\omega)$, if there exist positive numbers $\alpha_0$,
$\beta_0$ and an Riemann metric $g_0$ on $V$ such that
\begin{equation}
\omega(\xi,J_0\xi)\ge\alpha_0\|\xi\|^2_{g_0},\quad
|\omega(\xi,\eta)|\le\beta_0\|\xi\|_{g_0}\|\eta\|_{g_0},\quad
{\rm for}\;{\rm all}\;\xi,\eta\in TV;
\end{equation}
then there exists a positive number $\delta_0$ such that
$${\cal U}^m_{\delta_0}(J_0)=\bigl\{J_0{\rm exp}(-J_0Y)\bigm|
\,\|Y\|_{C^m}\le \delta_0,\;Y\in C^m({\rm T}_{J_0})\bigr\}
\subset{\cal J}^m_\tau(V,\omega)$$
for each integer $m\ge 1$. Here
$\|\cdot\|_{C^m}$ is defined in terms of the covariant derivatives
with respect to the Riemannian metric $g_0$. Furthermore,
$\delta_0>0$ can be chosen so small that every $J\in{\cal
U}^m_{\delta_0}(J_0)$ satisfies:
$\omega(\xi,J\xi)\ge\frac{\alpha_0}{2}\|\xi\|^2_{g_0}$ for all
$\xi\in TV$.}

\noindent{\it Proof}.\hspace{2mm}First note that the condition
(10) imply that
\begin{equation}
\frac{\alpha_0}{\beta_0}\|\xi\|_{g_0}\le\|J_0\xi\|_{g_0}\le\frac{\beta_0}{
\alpha_o}\|\xi\|_{g_0}
\end{equation}
for all $\xi\in TV$. Specially, we have that
$\alpha_0/\beta_0\le\|J_0(p)\|_{g_0}\le\beta_0/\alpha_0$ for all $p\in V$.

Next, for any $J=J_0{\rm exp}(-J_0Y)\in{\cal U}^m_{\delta_0}(J_0)$ and
$p\in V$, $\xi\in T_pV$ we have
\begin{eqnarray*}
\omega(\xi,J(p)\xi)&=&\omega(\xi, J_0(p){\rm exp}(-J_0(p)Y(p)))\\
&=&\omega(\xi,J_0(p)\xi)+\omega(\xi,J_0(p)[{\rm exp}(-J_0(p)Y(p))-I]\xi)\\
&=&\omega(\xi,J_0(p)\xi)+\omega(\xi,[{\rm exp}(J_0(p)Y(p))-I]J_0(p)\xi)\\
&\ge &\alpha_0\|\xi\|^2_{g_0}-\beta_0\|\xi\|_{g_0}
\|[{\rm exp}(J_0(p)Y(p))-I]J_0(p)\xi\|_{g_0}\\
&\ge &\alpha_0\|\xi\|^2_{g_0}-\frac{\beta_0^2}{\alpha_0}\|\xi\|^2_{g_0}
\|{\rm exp}(J_0(p)Y(p))-I\|_{g_0}.
\end{eqnarray*}
On the other hand, by (11) and the definition of {\bf exp}
\begin{eqnarray*}
\|{\rm exp}(J_0(p)Y(p))-I\|_{g_0}&\le &\sum^{\infty}_{k=1}
\frac{\|J_0(p)Y(p)\|^k_{g_0}}{k!}\\
 &\le &(\sum^{\infty}_{k=0}
\frac{\|J_0(p)Y(p)\|^k_{g_0}}{k!})\cdot{\rm exp}(\|J_0(p)Y(p)\|_{g_0})\\
&\le &\|J_0(p)\|_{g_0}\|Y(p)\|_{g_0}{\rm exp}(\|J_0(p)\|_{g_0}\|Y(p)\|_{g_0})\\
&\le &\frac{\beta_0}{\alpha_0}\delta_0\cdot{\rm exp}
(\frac{\beta_0}{\alpha_0}\delta_0).
\end{eqnarray*}
Thus we get
$$\omega(\xi, J(p)\xi)\ge [\alpha_0-\frac{\beta_0^3\delta_0}{\alpha_0}
{\rm exp}(\frac{\beta_0}{\alpha_0}\delta_0)]\|\xi\|^2_{g_0}.$$
Hence we can choose a positive number
$\delta_0\le\frac{1}{2}(\frac{\alpha_0}{\beta_0})^3$ so small that
$\alpha_0-\frac{\beta_0^3\delta_0}{\alpha_0} {\rm
exp}(\frac{\beta_0}{\alpha_0}\delta_0)\ge\alpha_0/2$. Lemma 2.7 is
proved. \hfill$\Box$\vspace{2mm}

Now  every ${\cal U}^m_{\delta_0}(J_0)$ is a Banach manifold, but it is not
separable or even has not a countable base. In order to be able to apply
Sard-Smale theorem in the transversality arguments below  we introduce
the space of the following type, which is one of our key techniques
in this paper.

Take a proper Morse function $h$ on $V$ and two sequences of regular values
of it, ${\rm a}=\{a_i\}$ and ${\rm b}=\{b_i\}$ satisfying:
$${\rm min}_{x\in V}h(x):= a_1<a_2<b_1<a_3<b_2<a_4<\cdots<a_k<b_{k-1}
<a_{k+1}<\cdots;$$
and denote $Q_i:=\{a_i\le h\le b_i\}$, $i=1,\cdots$, we have
\begin{equation}
V=\bigcup^{\infty}_{i=1}Q_i,\quad Q_i\cap Q_{i+2}=\emptyset,\quad
Int(Q_i)\cap Int(Q_{i+1})\ne\emptyset,\;i=1,\cdots.
\end{equation}
Moreover, every $Q_i$ is a smooth compact submanifold with smooth
boundary and has the same dimension as $V$. Following [F] we may
choose a sequence of sufficiently rapidly decreasing positive
numbers $\varepsilon^{(i)}=\{\varepsilon_k^{(i)}\}^{\infty}_{k=1}$
such that the space $C^\infty_{\varepsilon^{(i)}}({\rm
T}_{J_0}|_{Q_i})$ of those smooth sections $X\in C^\infty({\rm
T}_{J_0}|_{Q_i})$ for which
\begin{equation}
\|X\|_{\varepsilon^{(i)}}=\sum^\infty_{k=1}\varepsilon_k^{(i)}\|X\|_{
C^k(Q_i)}<\infty,
\end{equation}
is separable and dense in $L^1({\rm T}_{J_0}|_{Q_i})$.
In addition we always require that all $\varepsilon_1^{(i)}$ equal to $1$.
Let
\begin{equation}
C^\infty_{\varepsilon^{(i)}}({\rm T}_{J_0}^{(i)}):=\bigl\{X\in
C^\infty({\rm T}_{J_0})\bigm|\,{\rm supp}X\subset Q_i,\;
\|X\|_{\varepsilon^{(i)}}<\infty\bigr\}.
\end{equation}
This is a separable Banach space with respect to norm
$\|\cdot\|_{\varepsilon^{(i)}}$. We denote
\begin{equation}
{\cal L}_\varepsilon(J_0,h,{\rm a},{\rm b})
\end{equation}
by the space of all sequences ${\rm X}=(X_1,X_2,\cdots)$ with
$X_i\in C^\infty_{\varepsilon^{(i)}}({\rm T}_{J_0}^{(i)})$ and
\begin{equation}
\|{\rm X}\|_\varepsilon=\sum^\infty_{i=1}\|X_i\|_{\varepsilon^{(i)}}<
\infty.
\end{equation}
Then $({\cal L}_\varepsilon(J_0,h,{\rm a},{\rm
b}),\|\cdot\|_\varepsilon)$ is a separable Banach space. Let ${\rm
B}({\cal L}_\varepsilon(J_0,h,{\rm a},{\rm b});\delta)$ be a
closed ball in this Banach space of radius $\delta$. Then, for
sufficiently small $\delta\in (0, {\delta_0}/2)$ that $\|{\rm
X}\|_\varepsilon\le\delta$ implies that
$\|\sum^\infty_{i=1}X_i\|_{C^1}\le 2\delta$ and thus from Lemma
2.7 it follow that  $J:=J_0{\rm exp}(-J_0(\sum^\infty_{i=1}X_i))$
belongs to ${\cal J}_\tau(V,\omega)$ which is the space of all
smooth $\omega$-tame almost complex structures, and $(V,\omega, J,
g_0)$ is still g.bounded. Later, we fix such a $\delta$ and for
convenience denote by
$$\Xi:{\rm B}({\cal L}_\varepsilon(J_0,h,{\rm a},{\rm b});\delta)
\to{\cal J}_\tau(V,\omega),\;{\rm X}\mapsto J_0{\rm
exp}(-J_0(\sum^\infty_{i=1}X_i)).$$ and also  by
\begin{equation}
{\cal U}_\delta(J_0,h,{\rm a},{\rm b},\varepsilon)
\end{equation}
the image of
${\rm B}({\cal L}_\epsilon(J_0,h,{\rm a},{\rm b});\delta)$
under $\Xi$. This set is not necessary connected
in ${\cal J}_\tau(V,\omega)$.

Having the space many regularity results on compact symplectic
manifolds can be generalized to noncompact geometrically bounded
symplectic manifolds. In fact, the above construction can be
suitably modified so that the result of the moduli spaces in [Mc1]
may be generalized to any noncompact symplectic manifolds without
boundary, that is, the following proposition holds.

\noindent{\bf Proposition 2.8.}\hspace{2mm}{\it Given $A\in
H_2(V)$ and a closed Riemann surface $\Sigma$ of genus $g$ with
the complex structure $j$, and $J_0$ as above, then there is a
subset ${\rm B}_{\rm reg}({\cal L}_\varepsilon(J_0,h,{\rm a},{\rm
b}); \delta)$ of the second category in ${\rm B}({\cal
L}_\varepsilon(J_0,h,{\rm a},{\rm b}); \delta)$ such that for
every ${\rm X}\in{\rm B}_{\rm reg}({\cal L}_\varepsilon(J_0,h,{\rm
a},{\rm b}); \delta)$  the space
$${\cal M}_s(\Sigma, A, \Xi({\rm X}))$$
of all simple $\Xi({\rm X})$-holomorphic maps from $\Sigma$ to $V$
and representing $A$ is a smooth manifold of dimension $(1-g){\rm
dim}M+ 2c_1(A)$ and with a natural orientation. Moreover, for any
two ${\rm X}$ and ${\rm Y}$ in ${\rm B}_{\rm reg}({\cal
L}_\varepsilon(J_0,h,{\rm a},{\rm b}); \delta)$ it may be proved
that  ${\cal M}_s(\Sigma, A, \Xi({\rm X}))$ and
 ${\cal M}_s(\Sigma, A, \Xi({\rm Y}))$ are oriented cobordant.}

 For every integer $m>1$ we denote by
$\overline{Hom}^m_{J_0}(T\Sigma,TV)$ the space of the $C^m$-smooth sections
of the bundle of anti-$J_0$-linear homomorphisms from $T\Sigma$ to $TV$ over
$\Sigma\times V$. Consider the Banach vector bundle ${\cal H}^m$ over
${\rm B}({\cal L}_\varepsilon(J_0,h,{\rm a},{\rm b});\delta)$ whose
fibre at a point ${\rm X}$ is $\overline{Hom}^m_{\Xi({\rm X})}(T\Sigma,TV)$.
It is easy to know that this is a separable $C^m$ Banach vector bundle.
We call the elements of the bundle as inhomogeneous terms.
Fix a large integer $m_0>0$ such that the conditions of the Sard-Smale
theorem are satisfied. For every integer $m\ge m_0$ one may, as in [RT1],
prove that there exists a subset ${\cal H}^m_{\rm reg}$ of the second
category in ${\cal H}^m$ such that for every
$({\rm X},\nu)\in{\cal H}^m_{\rm reg}$ the space
\begin{equation}
{\cal M}^m_A(\Sigma,\Xi({\rm X}),\nu)
\end{equation}
of all $(\Xi({\rm X}),\nu)$-map from $\Sigma$ to $V$ representing $A$
is a $C^m$-smooth manifold of dimension $(1-g){\rm dim}M+ 2c_1(A)$ and with
a natural orientation. Moreover, for any two pairs $({\rm X},\nu)$ and
$({\rm Y},\mu)$ in ${\cal H}^m_{\rm reg}$
it may be proved that  ${\cal M}^m_A(\Sigma,\Xi({\rm X}),\nu)$ and
 ${\cal M}^m_A(\Sigma,\Xi({\rm Y}),\mu)$ are $C^m$
 oriented cobordant.  Specially, it should be noted that
\begin{equation}
{\cal H}^{m_0}_{\rm reg}\supseteq{\cal H}^{m_0+1}_{\rm reg}\supseteq\cdots,
\end{equation}
which implies that for any $({\rm X},\nu)\in{\cal H}^m_{\rm reg}$ and
$({\rm X}^\prime,\nu^\prime)\in{\cal H}^{m^\prime}_{\rm reg}$ with
$m^\prime>m\ge m_0$ the spaces ${\cal M}^m_A(\Sigma,\Xi({\rm X}),\nu)$ and
 ${\cal M}^{m^\prime}_A(\Sigma,\Xi({\rm X}^\prime),\nu^\prime)$ are also
 $C^m$ oriented  cobordant.

Let ${\cal G}(V)$ be the set of all Riemann  metrics on $V$ whose
injectivity radius are more than zero and sectional curvatures
have  upper bounds.
 We also denote by ${\cal GJ}_\tau(V,\omega)$ the set of all
$J\in{\cal J}_\tau(V,\omega)$ which satisfy:
$\omega(\xi,J\xi)\ge\alpha_0\|\xi\|^2_g$ and $|\omega(\xi,\eta)|\le\beta_0
\|\xi\|_g\|\xi\|_g$ for some fixed $g\in{\cal G}(V)$, constants
$\alpha_0$,$\beta_0>0$ and all $\xi,\eta\in TV$. Obviously,
for every g.bounded symplectic manifold $(V,\omega)$,
${\cal GJ}_\tau(V,\omega)$ is a nonempty  open subset of
${\cal J}_\tau(V,\omega)$ with respect to $C^1$-topology. However,
we do not affirm it to be connected.
For every connected component ${\cal G}(V)_c$ of ${\cal G}(V)$ we denote
${\cal GJ}_\tau(V,\omega)_c$ by the subset of $J\in{\cal GJ}_\tau(V,\omega)$
for which $(V,\omega,J,g)$ is g.bounded for some $g\in{\cal G}(V)_c$.
Using B\'evennec's construction [p.44,ALP] it is easily proved that
every ${\cal GJ}_\tau(V,\omega)_c$ is connected. Similarly, for every integer
$m\ge m_0$ we also denote by
\begin{equation}
 {\cal GH}^m(V,\omega):=\bigm\{\,(J,\nu)\bigm|\,
J\in{\cal GJ}_\tau(V,\omega)\;\; {\rm and}\;\;\nu\;\;{\rm satifying} (1)
\;\;\bigr\}
\end{equation}
and corresponding component ${\cal GH}^m(V,\omega)_c$, where
$\nu\in\overline{Hom}^m_J(TV, TM)$.
Then the later is still
connected.

\section{Transversality and Compactness}

In this section we shall follow the methods in [RT1][McSa1] to make arguments.
Because the techniques are same basically we only give the necessary
improvements and list the main results.

First of all, we start with the following notion. A pair
$(\Sigma;\bar{\bf z})$ of a connected Hausdorff topological space
$\Sigma$ and $k$ different points $\bar{\bf z}=\{z_1,\cdots,z_k\}$
on it is called the {\it semistable curve with $k$ marked
points}([FO]) if there exists a finite family of smooth closed
Riemann surfaces $\{\widetilde\Sigma_s:s\in\Lambda\}$ and
continuous maps
$\pi_{\scriptstyle\widetilde\Sigma_s}:\widetilde\Sigma_s\to\Sigma$
such that: $(i)$ each  $\pi_{\widetilde\Sigma_s}$ is a local
homeomorphism; $(ii)$ for each $p\in\Sigma$ it holds that
$1\le\sum_s \sharp\pi_{\widetilde\Sigma_s}^{-1}(p)\le 2$, and all
points which satisfy $\sum_s
\sharp\pi_{\widetilde\Sigma_s}^{-1}(p)=2$ are isolated; $(iii)$
for each $z_i$, $\sum_s
\sharp\pi_{\widetilde\Sigma_s}^{-1}(z_i)=1$. Denote by
$\Sigma_{\rm sing}:=\{p\in\Sigma: \sum_s
\sharp\pi_{\widetilde\Sigma_s}^{-1}(p)=2\}$ the set of all
singular points of $\Sigma$. Specially, each singular point $p$
such that $\sharp\pi_{\widetilde\Sigma_s}^{-1}(p)=2$ is called the
self-intersecting point of $\Sigma$. Call
$\Sigma_s:=\pi_{\widetilde\Sigma_s}(\widetilde\Sigma_s)$ the
$s$-th components of $\Sigma$, and $\widetilde\Sigma_s$ the smooth
resolution of $\Sigma_s$. Each $z_i$ is called the marked point.
The points in $\pi_{\widetilde\Sigma_s}^{-1}(\Sigma_{\rm sing})$
and $\pi_{\widetilde\Sigma_s}^{-1}(\bar{\bf z})$ are called the
singular points and  the marked points on $\widetilde\Sigma_s$,
respectively. Let $k_s$ be the number of all singular and marked
points on $\widetilde\Sigma_s$ and $g_s$ be the genus of
$\widetilde\Sigma_s$. The {\it genus} $g$ of $(\Sigma;\bar{\bf
z})$ is defined by
$$1 + \sum_s g_s + \sharp{\rm Inter}(\Sigma)-\sharp{\rm Comp}(\Sigma),$$
where $\sharp{\rm Inter}(\Sigma)$ and $\sharp{\rm Comp}(\Sigma)$
stand for the number of the intersecting points on $\Sigma$ and
the number of the components of $\Sigma$ respectively. If
$k_s+2g_s\ge 3$ we call the component
$(\widetilde\Sigma_s;\bar{\bf z}_s)$ stable. When all components
of $(\Sigma;\bar{\bf z})$ are stable we call $(\Sigma;\bar{\bf
z})$ the stable curve of genus $g$ and with $k$ marked points.

For the above genus $g$ stable curve $(\Sigma;\bar{\bf z})$ a
continuous map $f:\Sigma\to V$ is called $C^l$($l\ge 1$) if each
$f\circ\pi_{\widetilde\Sigma_s}$ is so. The homology class of $f$
is defined by
$f_*([\Sigma])=\sum_s(f\circ\pi_{\widetilde\Sigma_s})_*([\widetilde\Sigma_s])$.
An $C^m$ inhomogeneous term $\nu$ over $\Sigma$ is a set
$\{\nu_s:s\in\Lambda\}$ of inhomogeneous terms, where each $\nu_s$
is an $C^m$ inhomogeneous term of $\widetilde\Sigma_s$ and they
together satisfy the match conditions. A map $f:\Sigma\to V$ is
called $(J,\nu)$-perturbed holomorphic if each
$f\circ\pi_{\widetilde\Sigma_s}$  is $(J,\nu_s)$-perturbed
holomorphic. Denote by ${\cal M}^m_A(\Sigma,J,\nu)$ the moduli
space of all $(J,\nu)$-perturbed holomorphic maps from $\Sigma$
into $V$ with $f_*([\Sigma])=A$. Using the method in \S2 and the
arguments in [McSa1] [RT1.Prop.4.13] it follows that for every
given pair $(J^\prime,\nu^\prime)$
 with $C^m$ inhomogeneous term $\nu^\prime$ there exists
a pair $(J, \nu)$ with $C^m$ inhomogeneous term $\nu$
which may be arbitrarily close to it, such that  moduli space
${\cal M}^m_A(\Sigma,J,\nu)$ is a $C^m$-smooth manifold of dimension
$2c_1(V)(A)+2n(1-g)$.
In order to get suitable compactification of the above moduli space
the following form {\rm cusp-curve} due to Gromov was introduced in [RT1].
Given a $k$-point genus $g$ stable curve $(\Sigma;\bar{\bf z})$ as above,
$(\Sigma^\prime;\bar{\bf z}^\prime)$ is another $k$-point curve obtained
from it as follows:  First at some double points of
$\Sigma$ we join chains of $\mbox{\Bb C}P^1$ to separate the two components
and then attach some trees of $\mbox{\Bb C}P^1$, but require that
if one attaches a tree of $\mbox{\Bb C}P^1$ at a marked point
$x_i$, this $x_i$ will be replaced by a point different from intersection
points on some component of the tree, and under other cases the marked points
do not change. The components of $\Sigma$ is called {\it principal
components} and other {\rm bubble components}.
A continuous map  $f:\Sigma^\prime\to (V,\omega)$ is called a $\Sigma$-cusp
$(J,\nu)$-map if  for each principal component $\Sigma_s$
the map $f\circ\pi_{\widetilde\Sigma_s}$ is $(J,\nu_s)$-perturbed holomorphic
and the restriction of $f$ to a bubble component is a nonconstant
$J$-holomorphic map.
We define a $(\Sigma, J,\nu)$-{\rm cusp curve} as an equivalence class of cusp
maps modulo the parametrization groups of bubbles.  Its homology class is
defined as the sum of the homology classes of all components of the
any cusp map representatives of it. Denote by ${\cal CM}^m_A(\Sigma,J,\nu)$
the set of all $(\Sigma, J,\nu)$-cusp curves with the total homology class
$A$. For every element of the space one can obtain a reduced
$(\Sigma, J,\nu)$-cusp curve by forgetting multiplicity of the multiple
covering maps on bubble components and collapsing each subtree of the
bubbles whose components have the same image. Notice that this new cusp
curve may have different total homology class from the original one.
We denote by $\overline{\cal M}^m_A(\Sigma,J,\nu)$ the set of all reduced
$(\Sigma,J,\nu)$-cusp curve from ${\cal CM}^m_A(\Sigma,J,\nu)$.
For the semi-positive closed symplectic manifold $(V,\omega)$ it was proved
that  ${\cal CM}^m_A(\Sigma,J,\nu)$ is the cusp curve compactification
of ${\cal M}^m(\Sigma,J,\nu)$ and
$\overline{\cal M}^m_A(\Sigma,J,\nu)\setminus{\cal M}^m(\Sigma,J,\nu)$
consists of finitely many strata and each stratum is also branchedly covered
by a $C^m$-smooth manifold of codimension at least 2 ([RT1]). However,
in our case ${\cal CM}^m_A(\Sigma,J,\nu)$  is only the closure of
${\cal M}^m(\Sigma,J,\nu)$
due to the noncompactness of $(V,\omega)$. In order to get desire results
we assume that each $\nu_s$ in $\nu=\{\nu_s:s\in\Lambda\}$ satisfies (1).
For any compact subset $K\subset V$ let
${\cal CM}^m_A(\Sigma,J,\nu, K)$ be the subset of
${\cal CM}^m_A(\Sigma,J,\nu )$ consisting of all elements whose images are
intersecting with $K$.   Then we have

\noindent{\bf Proposition 3.1.}\hspace{2mm}{\it Let
$(V,\omega,g,J)$ be a g.bounded symplectic manifold and
$(\Sigma;\bar{\bf z})$ a $k$-point genus $g$ stable curve with a
bounded inhomogeneous term $\nu$ over it. Then there exists a
positive number $\eta=\eta(A,i(V,g),C_0,\alpha_0,\beta_0,\nu)$
such that
$$\bigcup_{f\in{\cal CM}^m_A(\Sigma,J,\nu,K)}{\rm Im}(f)\subset K_{\eta}.$$}

In fact, if $\Sigma_1,\cdots,\Sigma_p$ are principal components of
$\Sigma^\prime$ which only depend on $\Sigma$, and
$B_1,\cdots,B_l$ are bubble components of $\Sigma^\prime$ then it
follows from the proof of Lemma 4.5 in [RT1] that there is a
uniform constant $c$ such that $E(f_{\Sigma_i})\le c
(\omega(f_*([\Sigma_i]))+1)$, $E(f_{B_j})\le c \omega(f_*([B_j]))$
and therefore
$$\sum^p_{i=1}E(f_{\Sigma_i})+\sum^l_{j=1}b_j E(f_{B_j})\le c(\omega(A)+p).$$
Here the positive integer $b_1,\cdots,b_l$ satisfy
$A=\sum^p_{i=1}f_*([\Sigma_i])+ \sum^l_{j=1}b_j f_*([B_j])$. These
show that one can find a positive integer
$l_0=l_0(\omega(A),\Sigma,V,\omega, K)$ such that  it bounds $l$
uniformly. Moreover, for given area forms $\sigma_s$ on
$\widetilde\Sigma_s$($s=1,\cdots,p$) one can find a sufficiently
large $N>0$ such that all $(\widetilde\Sigma_s\times V,
\tilde\omega_s,\tau_s\oplus g,\tilde J_s)$ are g.bounded. Here
$\tilde\omega_s=N\sigma_s\times\omega$ and $\tilde J_s$ are
defined as in Proposition 2.6. From the proof of Lemma 2.5 it
follows that $\omega(f_*([\Sigma_s]))\ge
-N\int_{\Sigma_s}\sigma_s$ for each $s$. Combing these with
$\sum^p_{s=1}\omega(f_*([\Sigma_s]))\le\omega(A)$ we get that
$\omega(f_*([\Sigma_s]))\le \omega(A)+pN{\rm
min}_s\int_{\Sigma_s}\sigma_s$ and $\omega(f_*([B_j]))\le
p\omega(A)+ p(p-1)N{\rm min}_s\int_{\Sigma_s}\sigma_s$. Now since
$\Sigma^\prime$ is connected, by repeatedly using Lemma 2.5 we can
finish the proof of Proposition 3.1. As a consequence of this
proposition and Proposition 3.1 in [RT1] we have

\noindent{\bf Corollary 3.2.}\hspace{2mm}{\it For any compact
subset $K\subset V$, ${\cal CM}^m_A(\Sigma,J,\nu, K)$ is compact.}

As in [RT1,\S4],
$\overline{\cal M}^m_A(\Sigma,J,\nu)\setminus{\cal M}^m_A(\Sigma,J,\nu)$
can be stratified and their strata are indexed by
${\cal D}^{J,\nu}_{A,\Sigma}$(cf.[RT1] for definition). For a compact subset
$K\subset V$ we denote by ${\cal D}^{J,\nu}_{A,\Sigma}(K)$ the subset of
${\cal D}^{J,\nu}_{A,\Sigma}$ consisting of those
$D\subset{\cal D}^{J,\nu}_{A,\Sigma}$ which has a $\Sigma$-cusp
$(J,\nu)$-map representative intersecting with $K$.
Then carefully checking the proof of Lemma 4.5 in [RT1] we can prove

\noindent{\bf Lemma 3.3.}\hspace{2mm}{\it For any compact subsets
$K\subset V$, ${\cal D}^{J,\nu}_{A,\Sigma}(K)$ is a finite set.
But ${\cal D}^{J,\nu}_{A,\Sigma}$ may be a countable set.}

Corresponding to Theorem 4.2 in [RT1] we may use the argument method in
\S2 to get the following structure theorem.

\noindent{\bf Theorem 3.4.}\hspace{2mm}{\it Let $(V,\omega)$ be a
g.bounded semi-positive symplectic manifold, then there is a dense
subset ${\cal GH}^m_{\rm reg}(V,\omega)$  in ${\cal
GH}^m(V,\omega)$ such that for each pair $(J,\nu)\in{\cal
GH}^m_{\rm reg}(V,\omega)$, the complementary $\overline{\cal
M}^m_A(\Sigma,J,\nu)\setminus{\cal M}^m_A(\Sigma,J,\nu)$ consists
of at most countable many strata and each stratum is branchedly
covered by a $C^m$-smooth manifold of  codimension at least $2$.
Moreover, there are only finitely many strata of $\overline{\cal
M}^m_A(\Sigma,J,\nu)\setminus{\cal M}^m_A(\Sigma,J,\nu)$ which can
intersect with  $\overline{\cal M}^m_A(\Sigma,J,\nu,K)$ for every
compact subset $K\subset V$.}

More precisely, if for each $D\in{\cal D}^{J,\nu}_{A,\Sigma}$ we
denote by ${\cal M}^m_\Sigma(D,J,\nu)$ the space of all
$C^m$-smooth $(\Sigma,J,\nu)$-cusp curves such that the
homeomorphism type of its domain, homology class of each
component, components which have the same image are specified by
$D$. Then from Theorem 4.7 and Proposition 4.14 in [RT1] and the
arguments in \S2 we can obtain

\noindent{\bf Theorem 3.5.}\hspace{2mm}{\it For every $(J,\nu)$ in
a dense subset ${\cal GH}^m_{\rm reg}(V,\omega)$ of ${\cal
GH}^m(V,\omega)$ and a $D$ in ${\cal D}^{J,\nu}_{A,\Sigma}$ there
exists a $C^m$-smooth branched covering manifold ${\cal
N}^m_\Sigma(\bar D,J,\nu)$ of ${\cal M}^m_\Sigma(D,J,\nu)$ whose
dimension is not more than $2c_1(V)(A)+2n(1-g)-2k_D-2s_D$. Here
$k_D$ is the number of bubble components of $D$ and $s_D$ is the
number of marked points which are bubbling points. Moreover, for
any two pairs $(J,\nu)$ and $(J^\prime,\nu^\prime)$ in ${\cal
GH}^m_{\rm reg}(V,\omega)\cap{\cal GH}_{\rm reg}(V,\omega)_c$
there is a  path $(J_\tau,\nu_\tau)$ connecting $(J,\nu)$ and
$(J^\prime,\nu^\prime)$ in ${\cal GH}^m_{\rm reg}(V,\omega)_c$
such that $\cup_{t\in[0,1]}{\cal N}^m_\Sigma(\bar
D,J_\tau,\nu_\tau)\times\{t\}$ is a $C^m$-smooth cobordism.}

It should be noted that the manifolds  ${\cal M}^m_A(\Sigma,J,\nu)$ and
${\cal N}^m_\Sigma(\bar D,J,\nu)$ carry a canonical orientation.

Denote by  ${\cal BD}^{J,\nu}_{A,\Sigma}$ the subset of
${\cal D}^{J,\nu}_{A,\Sigma}$ whose elements contain the bubble
components. From Theorem 3.5 we have
\begin{equation}
\overline{\cal M}^m_A(\Sigma,J,\nu)\setminus{\cal M}^m_A(\Sigma,J,\nu)\subset
\bigcup_{D\in{\cal BD}^{J,\nu}_{A,\Sigma}}{\cal M}^m_\Sigma(D,J,\nu).
\end{equation}

\section{Gromov-Witten Invariants}

In this section we shall follow the method in [McSa1] to define the
Gromov-Witten invariants of Ruan-Tian's form---mixed invariants.
First of all, we recall some evaluation map. For a $k$-point genus $g$
stable curve $(\Sigma,\bar{\bf z})$, $\bar{\bf z}=(z_1,\cdots,z_k)$ and
integers $l> 0$ consider the $C^m$-smooth map
\begin{equation}
e^m_{(\Sigma,\bar{\bf z},J,\nu)}:{\cal M}^m_A(\Sigma,J,\nu)\times
\Sigma^l\mapsto V^k\times V^l=V^{k+l}
\end{equation}
 given by
$$(f;y_1,\cdots,y_l)\mapsto (f(z_1),\cdots,f(z_k);f(y_1),\cdots,f(y_l)).$$
For each $D\in{\cal BD}^{J,\nu}_{A,\Sigma}$ the similar  map
$$e^m_{(D,J,\nu)}:{\cal M}^m_\Sigma(D,J,\nu)\times (\Sigma^\prime)^l\mapsto
V^{k+l}$$ can be defined. For each
$D\in{\cal BD}^{J,\nu}_{A,\Sigma}$, let
$\pi^m_D:{\cal N}^m_\Sigma(\bar D,J,\nu)\to {\cal M}^m_\Sigma(D,J,\nu)$
be a branched covering defined below Definition 4.6 in [RT1]. The
composition maps $e^m_D=e^m_{(D,J,\nu)}\circ\pi^m_D$ satisfy
\begin{equation}
\bigcap_{S\subset{\cal M}^m_A(\Sigma,J,\nu)\times \Sigma^l
\,{\rm compact}}\overline{e^m_{(\Sigma,\bar{\bf z},J,\nu)}\bigl(
[{\cal M}^m_A(\Sigma,J,\nu)\times \Sigma^l]\setminus S\bigr)}
\subset\bigcup_{D\in{\cal BD}^{J,\nu}_{A,\Sigma}}{\rm Im}(e^m_D).
\end{equation}
These show that $e^m_{(\Sigma,\bar{\bf z},J,\nu)}$ is a
$C^m$-smooth pseudo-cycle. Let us recall the notion of the
pseudo-cycles introduced on the page 90 of [McSa1]. A
$k$-dimensional $C^m$-smooth {\it pseudo-cycle} in $V$ is a
$C^m$-smooth map $f:M\to V$ defined on an oriented $C^m$-smooth
$k$-dimensional manifold $M$(possibly noncompact) such that  the
boundary
$$f(M^\infty)=\bigcap_{S\subset M\,compact}\overline{f(M-S)}$$
of $f(M)$ is of dimension at most $k-2$, i.e.,there exists a
$C^m$-smooth manifold $W$ of dimension at most $k-2$ and a
$C^m$-smooth map $g:W\to V$ such that $f(M^\infty)\subset g(W)$.
If $\overline{f(M)}$ is also compact in $V$ then we call $f$ a
{\it strong pseudo-cycle}. Clearly, in a compact manifold these
two notions are equivalent. According to the definition  the
identity map $V\to V$ is not a strong pseudo-cycle in the
noncompact manifold $V$.
>From Remark 7.1.1 in [McSa1] it easily follows that every integral homology
class $\alpha\in H_2(V,\mbox{\Bb Z})$ can be represented by a
$C^\infty$ strong pseudo-cycle $f:M\to V$. Every strong
pseudo-cycle determines a homology class, and bordant
pseudo-cycles determine the same homology class. But in the
noncompact manifold $V$ a pseudo-cycle does not necessarily
determine a homology class as the identity map from $V$ to $V$.
Moreover, it is easily checked that the product of two (strong)
pseudo-cycles is also a (strong) pseudo-cycle in the product
manifold. If $f_k:M\to V_k$ are (strong) pseudo-cycles, $k=1,2$,
then the map $M\to V_1\times V_2,\,m\mapsto (f_1(m), f_2(m))$ is
also (strong) pseudo-cycle. In \S5 below we will need these
conclusions.
  Two pseudo-cycles $e:P\to V$ and $f:Q\to V$ are called
{\it transverse} if either $e(P)\cap f(Q)=\emptyset$ or
$e(P^\infty)\cap\overline{f(Q)}=\emptyset$, $\overline{e(P)}\cap
f(Q^\infty)=\emptyset$ and $T_xV={\rm Im}de(p)+{\rm Im}df(q)$
whenever $e(p)=f(q)=x$. However, it should be noted that for two
transverse pseudo-cycles $e$ and $f$ as above, if one of them is a
strong pseudo-cycle
 $\Delta(e,f):=\{(p,q)\in P\times Q|
e(p)=f(q)\}$ is  a {\it compact} manifold of dimension ${\rm
dim}P+{\rm dim}Q-{\rm dim}V$. This statement can be derived from
the definition of transversality of pseudo-cycles directly.
 Specially, it is a finite set if
$P$ and $Q$ are of complementary dimension. Under our case Lemma
7.1.2 in [McSa1] are not applicable due to the noncompactness of
the manifold $V$, which implies that ${\rm Diff}^r(V)$ is not
separable Banach manifold for every integer $r>0$. We must give
its suitable modification form. This can be obtained with our
method in \S2.

Fix a large integer $r>0$ and as in (12) we denote by
\begin{equation}
\chi^r_i:=\{X\in\chi^r(V)|\,{\rm supp}X\subset Q_i,\;\|X\|_{C^r}<\infty\},
\end{equation}
where $\chi^r(V)$ are the space of all $C^r$-vector fields on $V$, and
$$\|X\|_{C^r}={\rm sup}_{x\in V}|X(x)|_g+ {\rm sup}_{x\in V}|\nabla_gX(x)|_g
+\cdots+{\rm sup}_{x\in V}|\nabla^r_gX(x)|_g,$$
$\nabla_g$ is the Levi-Civita connection of metric $g$.
Then every $(\chi^r_i,\,\|X\|_{C^r})$ is separable Banach space.
Denote by
\begin{equation}
\chi^r(V)_0
\end{equation}
the space of all sequences ${\rm X}=(X_1,X_2,\cdots)$ with
$X_i\in\chi^r_i$ and $$\|{\rm
X}\|_{gr}=\sum^\infty_{k=1}\|X_k\|_{C^r}<\infty.$$ Then it is
easily proved that $(\chi^r(V)_0,\|\cdot\|_{gr})$ is a separable
Banach space. Note that every ${\rm X}\in\chi^r(V)_0$ determines a
bounded $C^r$-smooth vector field, denoted by $\rho_r({\rm
X})=\sum^\infty_{i=1}X_i$. Clearly, the image of $\rho_r$ contains
all smooth vector fields with compact support on $V$. But every
$C^r$-smooth bounded vector field on complete Riemann manifolds
can uniquely determine a one-parameter $C^r$-smooth diffeomorphism
group. Let us denote by $\{F_t(\rho({\rm X})):t\in\mbox{\Bb R}\}$
the group determined by $\rho_r({\rm X})$. Define
\begin{equation}
{\cal F}^r:\chi^r(V)_0\to{\rm Diff}^r(V),\;
{\rm X}\mapsto F_1(\rho_r({\rm X})).
\end{equation}
It is easily checked that ${\cal F}^r$ is a $C^r$-smooth map.
Corresponding to Lemma 7.1.2 in [McSa1] we have the following
lemma.

\noindent{\bf Lemma 4.1.}\hspace{2mm}{\it If a $C^p$-smooth
pseudo-cycle $e:P\to V$  and a $C^q$-smooth one $f:Q\to V$ satisfy
\begin{equation}
\dim P+\dim Q\ge\dim V
\end{equation}
then
\begin{description}
\item[(i)] for every sufficiently large integer $r>{\rm min}\{p, q\}$
there exists a set
$\chi^r(V,e,f)\subset\chi^r(V)_0$ of the second category
such that $e$ is transverse to ${\cal F}^r({\rm X})\circ f$ for all
${\rm X}\in\chi^r(V,e,f)$; these $\chi^r(V,e,f)$ also satisfy:
$$\chi^r(V,e,f)\supseteq\chi^{r+1}(V,e,f)\supseteq\cdots,$$
which implies that for any ${\rm X}\in\chi^r(V,e,f)$ and
${\rm Y}\in\chi^s(V,e,f)$ with $s>r$ it holds that
$$({\cal F}^r({\rm X})\circ f)\cdot e= ({\cal F}^s({\rm Y})\circ f)\cdot e$$
provided that the equality in (27) also holds and one of $f$ and
$e$ is a strong pseudo-cycle;
\item[(ii)] if the equality in (27) holds, $e$ and $f$ are transverse and one of them is a strong
pseudo-cycle, then $\Delta(e,f)$ is a finite
set and in this case we denote by $\nu(x, y)$ the intersection number of
$e$ and $f$ at $(x, y)\in\Delta(e,f)$, and define
$$e\cdot f=\sum_{(x, y)\in\Delta(e,f)}\nu(x, y);$$
\item[(iii)] the intersection number $e\cdot f$ depends only on the bordism
classes of $e$ and $f$ when one of them varies in the bordism class of
the strong pseudo-cycle.
\end{description}}

\noindent{\it Proof}.\hspace{2mm}The proof can be finished as in
[McSa1]. We only need prove that the map
\begin{equation}
\Theta:P\times Q\times\chi^r(V)_0\to V\times V:(p,q,{\rm X})\mapsto
\bigl(e(p), {\cal F}^r({\rm X})(f(q))\bigr)
\end{equation}
is transverse to the diagonal $\triangle_V$.
For any $(p,q,{\rm X})\in\Theta^{-1}(\triangle_V)$ the differential of
$\Theta$ at it is given by
$$D\Theta(p,q,{\rm X})(\xi,\eta,{\rm Y})=\biggl(De(p)(\xi),\;
D\bigl({\cal F}^r({\rm X})\circ f\bigr)(q)(\eta)+
\bigl[D{\cal F}^r({\rm X})({\rm Y})\bigr](f(q))\biggr),
$$
where $(\xi,\eta,{\rm Y})\in T_pP\times T_qQ\times\chi^r(V)_0$.
Let $m=e(p)={\cal F}^r({\rm X})(f(q))$. For any given
$(u,v)\in T_{(m,m)}(V\times V)$ we wish to find $w\in T_mV$ and
$(\xi,\eta,{\rm Y})\in T_pP\times T_qQ\times\chi^r(V)_0$ such that
\begin{eqnarray}
&De(p)(\xi)=w+u\\
&D\bigl({\cal F}^r({\rm X})\circ f\bigr)(q)(\eta)+
\bigl[D{\cal F}^r({\rm X})({\rm Y})\bigr](f(q))=w+v.
\end{eqnarray}
By taking $\xi=0$, $\eta=0$ and $w=-u$ we need only find ${\rm
Y}\in\chi^r(V)_0$ such that
\begin{equation}
\bigl[D{\cal F}^r({\rm X})({\rm Y})\bigr](f(q))=v-u.
\end{equation}
For $f(q)=[{\cal F}^r({\rm X})]^{-1}(m)$, by definition of
${\cal F}^r({\rm X})$,
it is $\alpha_{\rm X}(1)$, where $\alpha_{\rm X}(t)$ is the unique solution
of the initial value problem
\begin{equation}
\dot\alpha_{\rm X}(t)=\bigl(\sum^\infty_{i=1}X_i\bigr)(\alpha_{\rm X}(t)),
\quad\alpha_{\rm X}(0)=f(q).
\end{equation}
For $s\in(-1,1)$ and ${\rm Y}\in\chi^r(V)_0$ we denote by
$\alpha_{{\rm X}+s{\rm Y}}(t)$ the unique solution of the initial problem
\begin{equation}
\dot\alpha_{{\rm X}+s{\rm Y}}(t)=\bigl(\sum^\infty_{i=1}X_i+sY_i\bigr)
(\alpha_{{\rm X}+s{\rm Y}}(t)),\quad
\alpha_{{\rm X}+s{\rm Y}}(0)=f(q).
\end{equation}
Then we need to find ${\rm Y}\in\chi^r(V)_0$ such that
\begin{equation}
\frac{d}{ds}\alpha_{{\rm X}+s{\rm Y}}(1)|_{s=0}=v-u.
\end{equation}
By localization method it is easy to find a smooth vector field
$Z$ with compact support on $V$ such that for the unique solution
curve family $\beta(\rho_r({\rm X})+sZ)(t)$ of $\rho_r({\rm
X})+sZ$ with initial value $f(q)$ at zero it holds that
$$\frac{d}{ds}\beta(\rho_r({\rm X})+sZ)(1)|_{s=0}=v-u.$$
Now using unit the decomposition technique it is easy to find a
${\rm Y}\in\chi^r(V)_0$ with $Z=\rho_r({\rm Y})$. Thus we prove
the transversality.

Moreover, the standard computation shows that the restriction of the natural
projection $\Pi$ from $P\times Q\times\chi^r(V)_0$ to $\chi^r(V)_0$ to
$\Theta^{-1}(\triangle_V)$ is a Fredholm operator with index
\begin{equation}
Index(\Pi|\Theta^{-1}(\triangle_V))={\rm dim}P+{\rm dim}Q-{\rm dim}V,
\end{equation}
which is only dependent on dimension of $P$,  $Q$ and $V$. Under
our assumption this index is less than or equal to zero. Thus we
need only fix an integer $r>0$ such that Sard-Smale theorem can be
applied. The remainder of the arguments are the same as that in
[McSa1]. \hfill$\Box$\vspace{3mm}

 Now let a $k$-point genus $g$ stable curve $(\Sigma,\bar{\bf z})$,
$A\in H_2(V,\mbox{\Bb Z})$ and the pair $(J,\nu)$ satisfy the regularity
requirements  in \S2 and \S3.
The integral homology classes $\{\alpha_i\}_{1\le i\le k}$ and
$\{\beta_j\}_{1\le j\le l}$ of $V$ satisfy
\begin{equation}
\sum^k_1(2n-{\rm deg}(\alpha_i))+\sum^l_1(2n-{\rm deg}(\beta_j)-2)
=2c_1(V)(A)+2n(1-g).
\end{equation}
We choose strong pseudo-cycles $f_i:P_i:\to V$ and $h_j:Q_j\to V$ representing
$\alpha_i$ and $\beta_j$($1\le i\le k,1\le j\le l$), respectively. Then
\begin{equation}
f:=\prod^k_{i=1} f_i\times\prod^l_{j=1} h_j:\prod^k_{i=1}P_i\times
\prod^l_{j=1}Q_j\to V^{k+l}
\end{equation}
is a strong pseudo-cycle representing the integral homology class
$\prod_i\alpha_i\times\prod_j\beta_j\in H_*(V^{k+l},\mbox{\Bb
Z})$. Since the compositions $f\circ\phi$ of this $f$ with any
$\phi\in{\rm Diff}^r(V^{k+l})$ are also $C^r$-strong pseudo-cycles
representing the same class, using Lemma 4.1 we can assume that
$f$ is transverse to
 $e^m_{(\Sigma,\bar{\bf z},J,\nu)}$ and all $e^m_{(D,J,\nu)}$ because of
 the countability  of ${\cal BD}^{J,\nu}_{A,\Sigma}$. By Lemma 4.1 and (36)
 we can define   the mixed invariant
 \begin{equation}
 \Phi_{(A,\omega,g)}(\alpha_1,\cdots,\alpha_k|\beta_1,\cdots,\beta_l)=
 f\cdot e^m_{(\Sigma,\bar{\bf z},J,\nu)}.
\end{equation}
In the case that (36) does not hold we also define
\begin{equation}
\Phi_{(A,\omega,g)}(\alpha_1,\cdots,\alpha_k|\beta_1,\cdots,\beta_l)=0.
\end{equation}

As in [RT1] we can use the arguments in \S2 and \S3 to prove that
$\Phi_{(A,\omega,g)}(\alpha_1,\cdots,\alpha_k|\beta_1,\cdots,\beta_l)$
is independent of choices of $(J,\nu)$ in a dense subset of ${\cal
HJ}^m(V,\omega)_c$, marked points $z_1,\cdots,z_k$ in $\Sigma$,
the conformal structures on $\Sigma$, sufficiently large integers
$r, m$ and strong pseudo-cycles $(P_i,f_i)$, $(Q_j,h_j)$
representing $\alpha_i$, $\beta_j$ for a given component ${\cal
HJ}^m(V,\omega)_c$ of ${\cal HJ}^m(V,\omega)$. For two different
components we do not know what relationships there are between
corresponding invariants. When talking about some property of the
invariants we always mean them to be with respect to some fixed
component without special statements.
 Similarly, the corresponding results to
Proposition 2.4, 2.5 and 2.6, 2.7 in [RT1] can be proved. In
particular, under our assumptions one can define the invariant
$\Phi_{(A,\omega,\cal C)}$  as \S7 in [RT1] and prove the {\it
composition law}:
\begin{equation}
\Phi_{(A,\omega,g)}(\alpha_1,\cdots,\alpha_k|\beta_1,\cdots,\beta_l)=
\Phi_{(A,\omega,\cal C)}(\alpha_1,\cdots,\alpha_k|\beta_1,\cdots,\beta_l)
\end{equation}
where ${\cal C}=(\Sigma,\bar{\bf z})$ is a $k$-point genus $g$ stable curve and
$\alpha_1,\cdots,\alpha_k,\beta_1,\cdots,\beta_l$ are integral  homology
classes of $V$.

As to the deformation invariance of these invariants with respect
to the semi-positive deformation class of $\omega$ we introduce
the following notion of deformation equivalence. Two semi-positive
symplectic form $\omega_0$ and $\omega_1$ on a (noncompact)
geometrically bounded symplectic manifold $V$ is called {\it
deformedly equivalent} if there exists a smooth  1-parameter
family of semi-positive symplectic forms $\omega_t$ connecting
$\omega_0$ and $\omega_1$, and a family of almost complex
structures $J_t$ such that all $(V,\omega_t, J_t, g)$ are
uniformly geometrically bounded with respect to some metric
$g\in{\cal G}(V)$, that is, there exist constants $\alpha_0$ and
$\beta_0$ such that two inequalities in $1^\circ$ of
 Definition 2.3 hold uniformly for all $\omega_t$. As usual we may
 use the above method to prove our Gromov-Witten invariants are
 invariant under such semi-positive deformations of $\omega$.

\noindent{\bf Example 4.2.}\hspace{2mm}For any closed manifold $N$
and any closed $2$-form $\Omega$ on $N$ consider the symplectic
manifold $(M,\omega)=(T^\ast N, \omega_{\rm can}+\pi^\ast\Omega)$
then for any $k$-point genus $g$ stable curve $(\Sigma,\bar{\rm
z})$, $A\in H_2(V,\mbox{\Bb Z})$, the integral homology classes
$\{\alpha_i\}_{1\le i\le k}$ and $\{\beta_j\}_{1\le j\le l}$ of
$V$ we have
$$\Phi_{(A,\omega,g)}(\alpha_1,\cdots,\alpha_k|\beta_1,\cdots,\beta_l)=0.$$

In fact, take any Riemannian metric $h$ on $N$ and denote by $H$
the induced Riemannian metric on $T^\ast N$ by $h$. Then from
proof of Proposition 4.1 in [Lu2] it easily follows that all
symplectic manifolds $(T^\ast N, \omega_t)$ are uniformly
geometrically bounded with respect to $H$. Here
$\omega_t=\omega_{\rm can}+t\Omega$, $t\in [0, 1]$. Furthermore,
the proof there also shows that one can take a smooth family of
almost complex structures $J_t$ such that every $J_t$ is
$\omega_t$-compatible and $(M, \omega_t, J_t, H)$ are uniformly
geometrically bounded. Now Chern class $c_1(TM, J_t)$ is
independent of $t$ and thus they are all zero because $c_1(TM,
J_0)=0$ is clear. Hence the symplectic forms $\omega_0=\omega_{\rm
can}$ and $\omega_1$ are deformedly equivalent. But it is clear
that $\Phi_{(A,\omega_{\rm
can},g)}(\alpha_1,\cdots,\alpha_k|\beta_1,\cdots,\beta_l)$ always
vanishes. The above deformation invariance leads to the
conclusion.

In order to define the quantum homology
\footnote{P. Seidel pointed out that in the orginal version using the Poincar\'e
duality on noncompact manifolds does not give rise to a product on
$H^\ast_c(V)$, and should consider the quantum homology of $(V,\omega)$ instead.}
 of $(V,\omega)$ we need to assume that
 \begin{equation}
\Gamma_\omega=H^S_2(V,\mbox{\Bb Z})/H^S_2(V,\mbox{\Bb Z})_0\;{\rm
is}\;{\rm finitely}\; {\rm generated},\end{equation}
 so that the Novikov ring
$\Lambda_\omega$ associated to the homomorphism
$\omega:\Gamma_\omega\to\mbox{\Bb R}$ is well-defined. Here
$H^S_2(V,\mbox{\Bb Z})_0$ is the subgroup of classes $\alpha$ in
$H^S_2(V,\mbox{\Bb Z})$ such that $\langle
[\omega],\alpha\rangle=0$ and $\langle
c_1(V,\omega),\alpha\rangle=0$. As usual we denote by
$$QH_\ast(V)=H_\ast(V)\otimes\Lambda_\omega,$$
where $H_\ast(V)$ stands for the quotient of $H_\ast(V,\mbox{\Bb
Z})$ modulo torsion. The quantum intersect product is given by
$$\alpha*_V\beta=\sum_{A\in\Gamma_\omega}(\alpha*_V\beta)_A\otimes e^A\in QH_{k+l-2n}(V)$$
for $\alpha\in H_k(V)$ and $\beta\in H_l(V)$, where
$(\alpha*_V\beta)_A\in H_{k+l+2c_1(A)-2n}(V)$ is determined by
$$(\alpha*_V\beta)_A\cdot_V\gamma=\Phi_{(A,\omega,0)}(\alpha,\beta,\gamma)\;{\rm for}\;{\rm all}\;
\gamma\in H_\ast(V).$$
This gives an ring structure on $QH_\ast(V)$.

\noindent{\bf Remark 4.3.}\hspace{2mm} For given  integral
homology classes $\alpha_1,\cdots,\alpha_k,\beta_1,\cdots,\beta_l$
and their strong pseudo-cycles representatives $f_i: P_i\to V$,
$h_j: Q_j\to V$ as in (37) it follows from $V$ being noncompact
g.bounded that there exist the diffeomorphisms $\phi\in{\rm
Diff}^r(V)$ such that the images of $f$ and $\hat\phi\circ
f:=\prod^k_i\phi\circ f_i\times\prod^l_j\phi\circ h_j$ are not
intersecting each other and even have the larger distances. But
$\hat\phi\circ f$ and $f$ are representing the same homology
classes, therefore from our results that if their Gromov-Witten
invariants are not zero then the maps in ${\cal
M}_A(\Sigma,J,\nu)$ are distributed over $V$ in an even way. In
the same time this seems also to show the complexity of the
distributions of the holomorphic curves in the general noncompact
symplectic manifolds.

\section{Gromov-Witten Invariants of Compact Symplectic Manifolds with
Contact Type Boundary}

Let $(V,\omega)$ is a $2n$-dimensional compact symplectic manifold
with contact type boundary $\partial V$. That is, there is a
one-form $\alpha$ on $\partial V$ such that
$d\alpha=\omega|_{\partial V}$ and $\alpha\wedge(d\alpha)^{n-1}$
is a volume form on $\partial V$. Such a form $\alpha$ is called a
{\it contact form} on $\partial V$. One can associate a noncompact
symplectic manifold $(\widetilde V,\widetilde\omega)$ as follows:
$$\widetilde V=V\bigcup_{{\partial V}\times\{1\}}\partial V\times [1,+\infty),
\qquad \widetilde\omega=\cases{\omega &on $V$;\cr d(t\alpha) &on
$\partial V\times [1,+\infty)$.\cr}$$ Here $t$ is the second
coordinate. For a $J\in{\cal J}(V,\omega)$ and a Riemannian metric
$h$ on $V$ we may extend them to $\widetilde J$ and $\widetilde h$
respectively so that $\widetilde J$ and $\widetilde h$ are
constant on the $\partial V\times\{t\}$. It is easily checked that
$(\widetilde V,\widetilde\omega,\widetilde J,\widetilde h)$ is a
g.bounded symplectic manifold. Moreover the inclusion
$i:V\to\widetilde V$ induces clear isomorphisms
$i_\ast:H_\ast(V,\mbox{\Bb Z})\to H_\ast(\widetilde V,\mbox{\Bb
Z})$ and $i^\ast: H^\ast(\widetilde V,\mbox{\Bb Z})\to
H^\ast(V,\mbox{\Bb Z})$. It is clear that $i^\ast(c_1(\widetilde
V,\widetilde J))= c_1(V,J)$ and
$i^\ast([\widetilde\omega])=[\omega]$. Consequently, $(\widetilde
V,\widetilde\omega)$ is semi-positive if only and if $(V,\omega)$
is semi-positive. For a class $\alpha\in H_\ast(V,\mbox{\Bb Z})$
we denote $\widetilde\alpha$ by $i_\ast(\alpha)$. Then for a given
$k$-point genus $g$ stable curve $(\Sigma,\bar z)$, $A\in H_2(V)$
and integral homology classes $\{\alpha_i\}_{1\le i\le k}$ and
$\{\beta_j\}_{1\le j\le l}$ of $V$ satisfying (36) we define
$$\Phi_{(A,\omega,g)}(\alpha_1,\cdots,\alpha_k|\beta_1,\cdots,\beta_l):=
\Phi_{(\widetilde A,\widetilde\omega,
g)}(\widetilde\alpha_1,\cdots,
\widetilde\alpha_k|\widetilde\beta_1,\cdots,\widetilde\beta_l).\leqno(42)$$
Since both the space of all Riemannian metrics on $V$ and ${\cal
J}(V,\omega)$ are contractible it is easy to check that the left
of (42) is independent on the choices of $J$ in a dense subset of
${\cal J}(V,\omega)$, marked points $z_1,\cdots, z_k$ in $\Sigma$
and conformal structures on $\Sigma$. Moreover, they also satisfy
the axioms that Gromov-Witten invariants satisfy on closed
symplectic manifolds. Notice that $(\widetilde
V,\widetilde\omega)$ always satisfies the assumption in (41). One
may naturally define a quantum ring
$QH_\ast(V)=H_\ast(V)\otimes\Lambda_\omega$ from (42) and the
agruments above Remark 4.3.

\section{Rigidity of the Loops in the Group of Hamiltonian Diffeomorphisms with Compact Support}

The quantum homology had been used to study the topology of
symplectomorphism groups and Hamiltonian symplectomorphism groups
on closed symplectic manifolds in [Se1] [Le] [LMP]. In this
section we will use the techniques developed in the previous
sections and their idea to study  these groups on noncompact
g.bounded symplectic manifolds. Without special statements our
$2n$-dimensional symplectic manifold $(V,\omega)$ is always
assumed to satisfy the following condition:
$$A\in\pi_2(V),\; 2-n\le c_1(A)<0\Longrightarrow\omega(A)\le 0.\leqno(43)$$
Given an element $\phi\in \pi_1({\rm Diff}(V),id)$ and any
a loop $S^1\to{\rm Diff}(V),\,t\mapsto\phi_t$ representing it one can define
an endomorphism $\partial_\phi: H_\ast(V,\mbox{\Bb Q})\to H_{\ast +1}(V,\mbox{\Bb Q})$
by setting $\partial_\phi([C])$ as a homology class represented by
the cycle $S^1\times C\to V,\,(t,x)\mapsto \phi_t(x)$ for a cycle $C$ in $V$.
The main result in [LMP] is that for a loop $\phi$ in the group ${\rm Ham}(V,\omega)$
the endomorphism $\partial_\phi$ vanishes identically if a $2n$-dimensional closed symplectic
manifold $(V,\omega)$ satisfies (43).
In this section we generalize their result as follows:

\noindent{\bf Theorem 6.1.}\hspace{2mm}{\it If a $2n$-dimensional
g.bounded symplectic manifold $(V,\omega)$ satisfies (41) (43)
then for any loop $\phi$ in ${\rm Ham}^c(V,\omega)$\footnote{P.
Seidel had constructed an example with a nontrival Hamiltonian
loop with compact support.} the endomorphism $\partial_\phi$
vanishes.}

Let ${\it GS}(V)$ be the set of the symplectic structures $\omega$
on $V$ satisfying (41) (43). For any $\omega\in{\it GS}(V)$ we
denote by $S_\omega: \pi_1({\rm Symp}^c_0(V,\omega))\to\pi_1({\rm
Diff}(V), id)$ and $H_\omega: \pi_1({\rm
Ham}^c_0(V,\omega))\to\pi_1({\rm Diff}(V), id)$ the homomorphisms
induced by the group inclusions respectively. As in [LMP], as a
consequence of Theorem 6.1 we get the following result on the
rigidity of Hamiltonian loops.

\noindent{\bf Corollary 6.2.}\hspace{2mm}{\it For an element
$\phi$ in $\pi_1({\rm Diff}(V),id)$ if there exist $\omega_1$ and
$\omega_2$ in ${\it GS}(V)$ such that $\phi\in{\rm
Im}(H_{\omega_1})\cap{\rm Im}(S_{\omega_2})$ then it also belongs
to ${\rm Im}(H_{\omega_2})$.}

For a $2n$-dimensional compact smooth manifold
$M$ with nonempty boundary $\partial M$ we denote ${\rm Cont}(M)$ by
the set of all symplectic structures on it for which (43) holds and $\partial M$ is
of contact type. ${\rm Diff}(M,\partial M)$ denote the subgroup consisting
of all elements $F\in {\rm Diff}(M)$ whose restriction to $\partial M$ is the identity.
For a symplectic structure $\omega$ on $M$ we denote by the subgroups\vspace{-2mm}
$${\rm Symp}(M,\partial M,\omega):={\rm Diff}(M,\partial M)\cap{\rm Symp}(M,\omega),
\quad{\rm Ham}(M,\partial M,\omega):={\rm Diff}(M,\partial M)\cap{\rm Symp}(M,\omega).$$
By [Se2] these spaces may have infinitely many connected components.
Notice that  in Exercise 10.13 on the page 318 of [McSa2]
it was pointed out that for a noncompact symplectic manifold $(V,\omega)$
without boundary the flux homomorphism is still well-defined on
${\widetilde{\rm Symp}}^c_0(V,\omega)$ and the corresponding result to Theorem 10.12 also holds when
${\rm Symp}_0(V,\omega)$ is replaced by ${\rm Symp}^c_0(V,\omega)$.
In fact, carefully checking the proof Theorem 10.12 in [McSa2] one can get the stronger conclusion
that for the isotopy
$$[0,1]\to {\rm Symp}_0^c(V,\omega),\;t\mapsto\psi_t$$
with $\psi_0=id$ and ${\rm Flux}(\{\psi_t\})=0$ one actually make it to be isotopic with
fixed endpoints to a Hamiltonian isotopy $\{\phi_t\}$ such that the support does not increase.
That is, if a compact subset $K\subset V$ is such that ${\rm Supp}\psi_t\subseteq K$
for all $t\in [0, 1]$, then $\{\phi_t\}$ may be required to satisfy:
${\rm Supp}\phi_t\subseteq K,\;\forall t\in [0, 1]$.
Using this remark and Corollary 6.2 we may obtain

\noindent{\bf Corollary 6.3.}\hspace{2mm}{\it For a
$2n$-dimensional compact smooth manifold $M$ with nonempty
boundary $\partial M$ and $\phi\in\pi_1({\rm Diff}(M,\partial M),
id)$ if ${\rm Cont}(M)$ is nonempty then for any two $\omega_1$
and $\omega_2$ in ${\rm Cont}(M)$ it holds that $\phi\in{\rm
Im}(H_{\omega_1})\cap{\rm Im}(S_{\omega_2})$ if and only if
$\phi\in{\rm Im}(H_{\omega_2})\cap{\rm Im}(S_{\omega_1})$.}

In fact, let $(\widetilde M,\widetilde\omega_1)$ and $(\widetilde M,\widetilde\omega_2)$
be the symplectic expansion as made in \S5 they obviously satisfy (41) (43).
Moreover, every element of ${\rm Diff}_0(M,\partial M)$ can be extended into one
of ${\rm Diff}_0^c(\widetilde M)$ by the identity extension. Thus
${\rm Symp}_0(M,\partial M,\omega_i)$ and ${\rm Ham}_0(M,\partial M,\omega_i)$ may be viewed
as the subgroups of ${\rm Symp}_0^c(\widetilde M,\widetilde\omega_i)$ and
${\rm Ham}_0^c(\widetilde M,\widetilde\omega_i)$, $i=1,2$, respectively.
Now the conclusion may be derived from Corollary 6.2.

As pointed out in [LMP] their results may be generalized to
arbitrary closed symplectic manifolds with the methods developed
in [FO], [LT1, LT2], [R3], [Sie]. However, as done in the previous
sections it seem to be very hard to generalize our results to
arbitrary noncompact g.bounded symplectic manifolds with their
methods.

 It is well-known that there exists a one-to-one
correspondence between elements of $\pi_1({\rm Symp}(V,\omega))$
and isomorphism classes of symplectic fibre bundles over $S^2$
with fibre $(V,\omega)$(cf.[LPM][Se1]). For a given loop
$\phi_{t\in [0,1]}$ in ${\rm Symp}(V,\omega)$ the correspondent
symplectic fibre bundle $P_{\phi}\to S^2$ may be obtained as
follows: let $D^+$ and $D^-$ be two copies of the closed disk
$D^2$ of radius $1$ of the plane bounded by $S^1$, one can glue
the trivial fibre bundles $D^{\pm}\times (V,\omega)$ by a map
$\Phi:\partial D^+\times V\to \partial D^-\times V: (2\pi
t,x)\mapsto (-2\pi t,\phi_t(x))$. According to [Se1] a symplectic
fibre bundle with fibre $(V,\omega)$ on $S^2$ is a smooth fibre
bundle $\pi:E\to S^2$ together with a smooth family
$\Omega=(\Omega_b)_{b\in S^2}$ of symplectic forms on its fibres
satisfying locally trivial condition and the transition function
taking its value in the group ${\rm Symp}(V,\omega)$. He also call
a symplectic fibre bundle $(E,\Omega)\to S^2$ as Hamiltonian if
there is a closed two-form $\widetilde\Omega$ on $E$ such that
$\widetilde\Omega|E_b=\Omega_b$ for all $b\in S^2$. Later, we call
such a closed two-form $\widetilde\Omega$ on the Hamiltonian fibre
bundle as {\it Hamiltonian form}.

Furthermore, from proof of Proposition 10.17 on the page 320 of [McSa2]
one can prove that for every loop $S^1\to {\rm Ham}^c(V,\omega),\,t\mapsto\phi_t$
there is a smooth function $H_\phi: S^1\times V\to$ with compact support to
generate it.
Especially, there is an exact sequence of groups
$$0\to\pi_1({\rm Ham}^c(V,\omega))\to\pi_1({\rm Symp}^c_0(V,\omega))\mathrel{\mathop\to^{
{\rm Flux}}}H^1_c(V,\mbox{\Bb R}),$$
where ${\rm Flux}$ is the flux homomorphism.
Consequently, from the proof of Proposition 2.9 in [Se1] it follows that
for a loop $\phi$ in ${\rm Symp}^c(V,\omega)$ the symplectic fibre bundle
$P_\phi\to S^2$ is Hamiltonian if and only if
the loop $\phi$ may be homotopic to a Hamiltonian loop in
${\rm Symp}^c(V,\omega)$, that is, a loop in ${\rm Ham}^c(V,\omega)$.

As in [LMP] using the Wang exact sequence of pair $(P_\phi, S^2)$:
$$\cdots\to H_{q-1}(V,\mbox{\Bb Z})\mathrel{\mathop\to^{\partial_{\phi\ast}}}H_q(V,\mbox{\Bb Z})
\mathrel{\mathop\to^{i_\ast}}H_q(P_\phi,\mbox{\Bb Z})\to H_{q-2}(V, \mbox{\Bb Z})\to\cdots$$
the proof of Theorem 6.1 can be reduced to the following equivalent theorem.

\noindent{\bf Theorem 6.4.}\hspace{2mm}{\it Let $(V,\omega)$ be as
in Theorem 6.1 and $\phi$ a loop in ${\rm Ham}^c(V,\omega)$. Then
the homomorphism $i:H_\ast(V,\mbox{\Bb Q})\to
H_\ast(P_\phi,\mbox{\Bb Q})$ is injective.}

In order to prove this theorem we need to give the detailed construction
in Proposition 2.9 of [Se1] since the more conclusions are needed.
 Let $D^+_{1/3}=\{z\in D^+|\,1/3\le |z|\le 1\}$
and $D^-_{1/3}=\{z\in D^-|\,1/3\le |z|\le 1\}$. Denote by $(r,t)_\pm$ the
polar coordinate in $D^\pm$ with $t\in S^1=\mbox{\Bb R}/\mbox{\Bb Z}$.
In the set $\Delta:=\{(r,t)_+, (r,t)_-\,|\,(r,t)\in D\}$
we define an equivalence relation $\sim$ as follows: the equivalence class
of $(r,t)_+$ is $[r,t]_+=\{(r,t)_+, (-r+5/3,-t)_-\}=[-r+5/3,-t]_-$ if $2/3\le r\le 1$,
those of $(r,t)_+$ and $(r,t)_-$ are $[r,t]_+$ and $[r,t]_-$ respectively if $0\le r<2/3$.
Then $S^2=\Delta/\sim$ and $U_\pm:=\{[r,t]_\pm\,|\,(r,t)\in D\}$ form an
open cover of $S^2$. $U_+\cap U_-=\{[r,t]_+=[-r+5/3,-t]_-\,|\,(r,t)\in [2/3,1]\times S^1\}$.
The coordinate charts $\varphi_\pm:D\to U_\pm,\,(r,t)\mapsto [r,t]_\pm$
give an atlas on $S^2$. The transition map is:
$$\varphi_-^{-1}\circ\varphi_+:
D_{1/3}:=\{z\in D\,|\,2/3\le |z|\le 1\}\to D_{1/3},\,(r,t)\to (-r+5/3,-t).$$
We also consider the formal set
$$\{((r,t)_\pm, x)\,|\,(r,t,x)\in D\times V\}\leqno(44)$$
and in it we define an equivalence relation
$\mathrel{\mathop\sim^{{\phi}}}$ as follows: the equivalence class
of $((r,t)_+,x)$ is $[r,t,x]_+^\phi=\{((r,t)_+,x)\}$ if $0\le
r<2/3$, that of $((r,t)_-, x)$ is $[r,t,x]_-^\phi=\{((r,t)_-,
x)\}$ if $0\le r<2/3$, and that of $((r,t)_+, x)$ is
$[r,t,x]_+^\phi=[-r+5/3, -t, \phi_t(x)]_-^\phi :=\{((r,t)_+, x),
((-r+5/3, -t)_-, \phi_t(x))\}$ if $2/3\le r\le 1$. Then the set,
denoted by $P_\phi$, of all equivalence classes of elements in the
set of (44) is a symplectic fibre bundle with fibre $(V,\omega)$.
Two  bundle charts $\Phi_+:U_+\times V\to
P_\phi|U_+,\,([r,t]_+,x)\mapsto [r,t,x]_+^\phi$ and
$\Phi_-:U_-\times V\to P_\phi|U_-,\,([r,t]_-,x)\mapsto
[r,t,x]_-^\phi$ form an bundle atlas of $P_\phi$. The transition
map is given by
$$\Phi_-^{-1}\circ\Phi_+: U_+\cap U_-\times V\to U_+\cap U_-\times V,
\,([r,t]_+,x)\mapsto ([-r+5/3,-t]_-, \phi_t(x)).\leqno(45)$$
Denote by $p_\pm:U_\pm\times V\to V$ the natural projections to
the second factors, and $\omega^\pm:=p_\pm^\ast\omega$. Define a
one-form $\theta_\phi$ on $U_+\times V$ as follows:
$\theta_\phi(([r,t]_+,x))=-\delta(r)H_\phi(t,\phi_t(x))dt$. Here
$H_\phi:S^1\times V\to\mbox{\Bb R}$ is a smooth function
generating $\phi_{t\in [0,1]}$ and having compact support,
$\delta:[0, 1]\to [0,1]$ is a monotone smooth function such that
$\delta(r)=0$ for $0\le r\le 1/4$ and $\delta(r)=r$ for $1/3\le
r\le 1$. In this paper we always fix this $\delta$ function.
 Clearly, $\theta_\phi$ has compact support.
Straightforward computation shows that the closed two-forms
$(\Phi_+^{-1})^\ast(\omega^+ + d\theta_\phi)$ and $(\Phi_-^{-1})^\ast \omega^-$
are the same on overlap $P_\phi|_{U_+\cap U_-}$.
Thus they define a closed two-form $\widetilde\Omega_\phi$ on $P_\phi$ by
$$\widetilde\Omega_\phi|_{P_\phi|U_+}=(\Phi_+^{-1})^\ast(\omega^+ + d\theta_\phi)\quad
{\rm
and}\quad\widetilde\Omega_\phi|_{P_\phi|U_-}=(\Phi_-^{-1})^\ast\omega^-.\leqno(46)$$
Let a compact subset $K\subseteq V$ be such that ${\rm
Supp}H_\phi\subseteq S^1\times K$. Then from the above definition
it easily follows that
$$P_\phi\setminus\Bigl(\Phi_+(U_+\times K)\bigcup\Phi_-(U_-\times K)\Bigr)=
S^2\times (V\setminus K),\leqno(47)$$ and on the set of (47) it
holds that
$$\widetilde\Omega_\phi=p_2^\ast\omega,\leqno(48)$$
where $p_2:S^2\times V\to V$ is the natural projection.
Moreover, one can easily prove that the above two-form $\widetilde\Omega_\phi$
is a Hamiltonian form on $P_\phi$ and
also satisfies:
$$\pi_\ast\widetilde\Omega_\phi^{n+1}=0\quad{\rm on}\;S^2\setminus\{[r,t]_+\in S^2\,|\,
1/4<r<1/3\},\leqno(49)$$ where $\pi_\ast$ is the fiber integration
map, and $\Omega_\phi$ a smooth family of symplectic forms on the
fibres of $P_\phi\to S^2$. Different from the case that $V$ is the
closed symplectic manifold we neither know the existence of a
Hamiltonian form $\widetilde\Omega$ on $P_\phi$ such that
$\pi_\ast\widetilde\Omega^{n+1}=0$ nor the uniqueness of such
forms if they exist. A Hamiltonian form $\widetilde\Omega$ on
$P_\phi$ is called to have {\rm CS} {\it property} if there are
compact subsets $K_\phi\subset V$ and $\widehat K_\phi\subset
P_\phi$ such that $P_\phi\setminus\widehat K_\phi=S^2\times
(V\setminus K_\phi)$ and on them it holds that
$\widetilde\Omega=p_2^\ast\omega$. Let us denote by
$${\cal H}(\phi)={\cal H}(P_\phi)$$
the set of all Hamiltonian forms having {\rm CS} property on
$P_\phi$. Since for any two Hamiltonian fibre bundles $P_\phi$ and
$P_\psi$ on $S^2$ obtained from loops $\phi_{t\in [0,1]}$ and
$\psi_{t\in [0, 1]}$ in ${\rm Ham}^c_0(V,\omega)$ one can always
find compact subsets $K\subset V$, $\widehat K_\phi\subset P_\phi$
and $\widehat K_\psi\subset P_\psi$ such that
$$P_\phi\setminus\widehat K_\phi= P_\psi\setminus\widehat K_\psi=S^2\times (V\setminus K),\leqno(50)$$
we may say a symplectic fibre bundle isomorphism ${\rm I}^{\phi\psi}$ between $P_\phi$ and $P_\psi$
to have {\rm CS} {\it property} if it is the identity map on the sets in (50),
that is, ${\rm I}^{\phi\psi}(z, v)=(z, v)$ for all $(z,v)\in S^2\times (V\setminus K)$.
Clearly, such an isomorphism induces a natural bijection
${\rm I}^{\phi\psi\ast}$ from ${\cal H}(\psi)$
to ${\cal H}(\phi)$ by the pull-back map.

For every $\widetilde\Omega\in {\cal H}(\phi)$ and the standard symplectic form $\sigma$ on $S^2$
it is easily proved that there is always a large constant $c(\widetilde\Omega, \phi)>0$
such that all two-forms $\widetilde\Omega + c\pi^\ast\sigma$
are symplectic forms on $P_\phi$ for all $c\ge c_\phi$.
Though these symplectic forms are also the Hamiltonian form on $P_\psi$,
but they have no {\rm CS} property.

Given a Hamiltonian form $\widetilde\Omega$ on $P_\phi$,
in [Se1] two continuous sections
$s_0$ and $s_1$ of $P_\phi$ are called $\Gamma_\omega$-equivalent if
$\widetilde\Omega_\phi(s_0)=\widetilde\Omega_\phi(s_1)$
and $c_1(TP_{\phi}^{\rm vert})(s_0)=c_1(TP_{\phi}^{\rm vert})(s_1)$.
The key point is this definition being independent of the choice of
$\widetilde\Omega$(cf.[Se1]).

 Following [Se1] we denote by ${\cal J}(P_\phi,\Omega_\phi)$ the space of smooth families
 ${\bf J}=(J_z)_{z\in S^2}$ of almost complex structures on the fibre of $P_\phi$
 such that $J_z$ is $\Omega_{\phi z}$-compatible for all $z$. In other words, ${\bf J}$
 is a smooth section of a bundle over $S^2$ whose fibre at a point $z\in S^2$
 is the space ${\cal J}(P_{\phi z}, \Omega_{\phi z})$. For the positively
 oriented complex structure $j$ on $S^2$ and ${\bf J}\in{\cal J}(P_\phi,\Omega_\phi)$,
 $\widehat{\cal J}(j,{\bf J})$ denote  the space of all almost complex structures
 $\hat J$ on $P_\phi$ compatible with $j$ and ${\bf J}$, that is, $\hat J$ satisfying:
 $D\pi\circ\hat J=j\circ D\pi$ and $\hat J|P_{\phi z}=J_z$ for all $z\in S^2$.
 Similarly, for every integer $m\ge 1$ we denote  $\widehat{\cal J}^m(j,{\bf J})$
  by the space of all $C^m$-smooth almost complex structures on $P_\phi$ satisfying
  the above conditions.
 A smooth section $s:S^2\to P_\phi$ is called $(j,\hat J)$-holomorphic if
 $ds\circ j=\hat J\circ ds$. For a given $\widetilde\Omega\in{\cal H}(\phi)$,
 from the above arguments  it is not difficult to
 prove that all symplectic manifolds $(P_\phi, \widetilde\Omega + c\pi^\ast\sigma)$
 are g.bounded with respect to some $\hat J\in\widehat{\cal J}(j,{\bf J})$ and
 some Riemannian metric on $P_\phi$. To see this point we choose a $g\in{\cal G}(V)$.
 Let $\tau_0$ be the standard metric on $S^2$. Notice that the above arguments show
 that one can choose a Riemannian metric  $G$ on $P_\phi$ such that it equals to
 $\tau_0\oplus g$ outside a compact subset.
 When $g$ takes over a connected component ${\cal G}(V)_c$ of ${\cal G}(V)$
 all corresponding Riemannian metrics on $P_\phi$ also form a connected subset of
 all Riemannian metrics on $P_\phi$, denoted by ${\cal G}(P_\phi)_c$. Later we always
 fix a component without special statements. For a $G\in {\cal G}(P_\phi)_c$
 we denote $G_z$ by the induced  metric on fibre $P_{\phi z}$ then one can use the
 standard method to find  ${\bf J}\in{\cal J}(P_\phi,\Omega_\phi)$ such that the
 family of symplectic manifolds  $\{(P_{\phi z}, \Omega_{\phi z}, G_z)\}_{z\in S^2}$
 is uniformly g.bounded. That is,
 their sectional curvature has a uniform upper bound, the injectivity radius has
 a uniform  positive lower bound and there exist positive constants $\alpha_0$ and
 $\beta_0$ such that
 $$\Omega_{\phi z}(\xi, J_z\xi)\ge\alpha_0\|\xi\|^2_{G_z}\qquad{\rm and}\qquad
 |\Omega_{\phi z}(\xi,\eta)|\le \beta_0\|\xi\|_{G_z}\|\eta\|_{G_z},\;\forall
 z\in S^2,\;\xi,\,\eta\in TP_{\phi z}.$$
 We denote by ${\cal GJ}(P_\phi,\Omega_\phi)_c$ all such ${\bf  J}\in{\cal J}(P_\phi,\Omega_\phi)$
 constructed from $\Omega_\phi$ and some $G\in{\cal G}(P_\phi)_c$ with the standard method.
 On the other hand from $(\widetilde\Omega + c\pi^\ast\sigma)|P_{\phi z}=\Omega_{\phi z}$
 and $G|P_{\phi z}=G_z$ it follows that the almost complex structure $\hat J$ on $P_\phi$
 constructed from $G$ and  $\widetilde\Omega + c\pi^\ast\sigma$ with the standard method
 must be in $\widehat{\cal J}(j,{\bf J})$ and such that
 $(P_\phi, \widetilde\Omega+ c\pi^\ast\sigma, \hat J, G)$
 is also g.bounded.
 Now fix such a ${\bf J}\in {\cal J}(P_\phi, \Omega_\phi)$ and a $\hat J\in\widehat{\cal J}(j,{\bf J})$,
  and as in \S2 we can construct a separable Banach space
 so that the transversity arguments in \S7 of [Se1] can be completed in our case.
 That is, under our assumptions, we can find $\hat J\in\widehat{\cal J}(j,{\bf J})$ such that
 \begin{description}
 \item[(i)]  $(P_\phi, \widetilde\Omega+ c\pi^\ast\sigma, \hat J, G)$ is g.bounded,
 \item[(ii)] the space  ${\cal S}(P_\phi, \Omega_\phi, j,\hat J, D)$ of all $(j,\hat J)$-holomorphic
 sections of $P_\phi$ representing a $\Gamma_\omega$-equivalence
class $D$ of a section of $P_\phi$ is a smooth manifold of dimension $2n+ 2c_1(TP_\phi^{\rm vert})(D)$
and for  chosen two different points $z_1,\, z_2\in S^2$ in advance and isomorphisms
$F_k^\phi:(P_{\phi z_k},\Omega_{\phi z_k})\to (V,\omega)$, $k=1, 2$, the maps
$${\rm EV}_k^{\phi D}:{\cal S}(P_\phi, \Omega_\phi, j,\hat J, D)\to V,\quad s\mapsto F_k^\phi(s(z_k))$$
are pseudo-cycles in the sense of \S7.1 of [McSa1].
\end{description}
Later we will fix such a $\hat J$ and a $c\ge c_\phi$ without special statements.
For two integral homology classes $\alpha,\beta\in H_\ast(V,\mbox{\Bb Z})$
and their strong pseduo-cycles representatives
$f_M:M\to V$ and $f_N:N\to V$ we can, as in \S4, show that there exist
$H\in{\rm Diff}(V\times V)$ such that the pseduo-cycle
${\rm EV}^{\phi D}:=({\rm EV}_1^{\phi D}, {\rm EV}_2^{\phi D})$ and strong pseudo-cycle $H\circ (f_M\times f_N)$
transversely intersect provided that
$$2n+ 2c_1(TP_\phi^{\rm vert})(D)+ {\rm deg}(\alpha)+{\rm deg}(\beta)=4n.\leqno(51)$$
Thus we may define a kind of Gromov-Witten invariants
$$\Phi_{(\phi,D; {\bf J})}(\alpha,\beta):={\rm EV}^{\phi D}\cdot (H\circ (f_M\times f_N))\leqno(52)$$
if (51) holds, and zero if (51) does not hold.
It is easy to prove that the right side of (52) is independent of
the choices of $\hat J$, $g$ ,$z_k$  and generic representatives.
When $\Gamma_\omega$ is finitely generated the rational
Novikov ring of it is well-defined
and thus quantum homology $QH_\ast(V)$ can be defined as in \S4.
In this case we use the idea from [LMP] to define the formal sum
$$\Psi^{\bf J}_{\phi,D}(\alpha)=\sum_{B\in\Gamma_\omega}\alpha_B\otimes e^B\leqno(53)$$
for $\alpha\in H_\ast(V,\mbox{\Bb Z})$, where $\alpha_B\in H_{\ast+d+2c_1(B)}(V)$ is determined by
$$\alpha_B\cdot_V\beta=\Phi_{(\phi,D+B; {\bf J})}(\alpha,\beta)\leqno(54)$$
for every $\beta\in H_\ast(V,\mbox{\Bb Z})$ and $B\in\Gamma_\omega$.
Here
$$d=c_1(TP_\phi^{\rm vert})(D)\leqno(55)$$
and $D+B$ is understood as in Lemma 2.10 of [Se1], that is, $D+B$ is
the only $\Gamma_\omega$-equivalence class of sections of $P_\phi$ such that
$$\widetilde\Omega_\phi(D+B)=\widetilde\Omega_\phi(D)+\omega(B)\qquad{\rm and}\qquad
c_1(TP_\phi^{\rm vert})(D+B)= c_1(TP_\phi^{\rm vert})(D) + c_1(B).$$
The following lemma shows that for every $\alpha\in H_\ast(V,\mbox{\Bb Z})$,
$\Psi_{\phi,D}(\alpha)$ is an element of $QH_{\ast+d}(V)$.

\noindent{\bf Lemma 6.5.}\hspace{2mm}{\it If $\Gamma_\omega$ is
finitely generated then for any $\alpha\in H_\ast(V,\mbox{\Bb
Z})$, $\Gamma_\omega$-equivalence class $D$ of sections of
$P_\phi$ and constant $C>0$ there are only finitely many $B\in
\Gamma_\omega$ such that $\alpha_B\ne 0$ and $\omega(B)\le C$ in
(53).}

\noindent{\it Proof}.\hspace{2mm} Since $\Gamma_\omega$ is
finitely generated the rational Novikov ring of it is a countable
set. Moreover, $\alpha_B\cdot_V\beta=0$ unless ${\rm
deg}(\alpha)+{\rm deg}(\beta)+ 2c_1(TP_\phi^{\rm vert})(D+B)=2n$.
Assume that there are a constant $C>0$ and infinitely many $B_i\in
\Gamma_\omega$ such that
$$\alpha_{B_i}\ne 0\quad{\rm and}\quad \omega(B_i)\le C,\; i=1,2,\cdots.$$
Then there are infinitely many $\beta_i\in H_\ast(V,\mbox{\Bb Z})$
such that
$$\Phi_{(\phi,D+B_i; {\bf J})}(\alpha,\beta_i)\ne 0,\qquad
{\rm deg}(\alpha)+{\rm deg}(\beta_i)+ 2c_1(TP_\phi^{\rm
vert})(D+B_i)=2n\leqno(56)$$ for all $i=1, 2,\cdots$. Recall the
definition in (52) we can always find $F_i\in{\rm Diff}(V\times
V)$ such that the image sets of all $F_i\circ (f\times h)$ are
contained in a fixed compact subset $S$ of $V\times V$. In fact,
from the proof of Lemma 4.1 one can find ${\rm
X}_i\in\chi^r(V\times V,{\rm EV}, f\times h)$ with $\|{\rm
X}_i\|_{gr}\le 1$ such that ${\rm EV}$ is transverse to all
$F_i\circ (f\times h)$ with $F_i:={\cal F}^r({\rm X}_i)$,
$i=1,2,\cdots$. But that $\|{\rm X}_i\|_{gr}\le 1$ implies that
$\|\rho_r({\rm X}_i)\|_{C^1}\le 2$ for all $i\ge 1$. Thus the
image sets of all maps $F_i\circ(f\times h)$ are contained in
 a fixed compact subset of $V\times V$, denoted by $S$.
The first formula of (56) shows that there exist $\hat J$-holomorphic section
$s_i$ representing the classes $D+ B_i$ with ${\rm EV}^{\phi(D+B_i)}(s_i)\cap S\ne\emptyset$.
In particular, there exists a compact subset $K$ of $P_\phi$ such that
$s_i(S^2)\cap K\ne\emptyset$ for all $i=1,2,\cdots$. Now
$$0\le (\widetilde\Omega+c\pi^\ast\sigma)(s_i)=\widetilde\Omega(D)+ \omega(B_i)
 +c\int_{S^2}\sigma\le\widetilde\Omega(D)+ c\int_{S^2}\sigma + C \leqno(57)$$
  because $D+B_i$ is the equivalence classes of sections of $P_\phi$ and
 $\int_{S^2}s^\ast(\pi^\ast\sigma)=\int_{S^2}(\pi\circ s)^\ast\sigma=\int_{S^2}\sigma$
 for every smooth section $s$ of $P_\phi$. This shows that
 there are infinitely many homology classes in $P_\phi$ with nonconstant
 $\hat J$-holomorphic spheres representatives whose image intersects with
 a fixed compact subset $S$ in $P_\phi$. It contradicts to Gromov compactness theorem.
 \hfill$\Box$\vspace{2mm}

 Consequently, (53) defines a $\Lambda_\omega$-linear homomorphism
 $\Psi^{\bf J}_{\phi,D}$ from $QH_\ast(M)$ to $QH_{\ast + d}(M)$ with
 $d=2c_1(TP_\phi^{\rm vert})(D)$.
 Moreover, if loops $\phi_{t\in [0,1]}$ and $\chi_{t\in [0,1]}$ are homotopic
 in ${\rm Ham}^c_0(V,\omega)$ there exists a Hamiltonian fibre bundle
 isomorphism ${\rm I}^{\phi\chi}$ having {\rm CS} property from
 $P_\phi$ to $P_\chi$. For a ${\bf J}\in{\cal J}(P_\phi,\Omega_\phi)$ and a $\Gamma_\omega$-
 equivalence class $D$ of sections of $P_\phi$ the isomorphism
 ${\rm I}^{\phi\chi}$ determines a ${\rm I}^{\phi\chi}_\ast({\bf J})$ and a
  $\Gamma_\omega$-equivalence class ${\rm I}^{\phi\chi}_\ast(D)$ of sections of $P_\chi$.
  It is not hard to prove that
  $$\Psi^{\bf J}_{\phi, D}=\Psi^{{\rm I}^{\phi\chi}_\ast({\bf J})}_{\chi, {\rm I}^{\phi\chi}_\ast(D)}.\leqno(58)$$
  As in [Se1][LMP] we have
 $$\Psi^{\bf J}_{\phi, D+B}=\Psi^{\bf J}_{\phi,D}\otimes e^{-B}$$
 for every $B\in\Gamma_\omega$ and the $\Gamma_\omega$-equivalence classes $D$ of sections
 of $P_\phi$, and the following conclusion.

 \noindent{\bf Lemma 6.6.}\hspace{2mm}{\it For the constant loop $\phi_0=\{id\}$ and
 the $\Gamma_\omega$ equivalence class $D_0$ of the flat section $s_0=S^2\times\{pt\}$
 of $P_{\phi_0}=S^2\times V$ the map $\Psi^{\bf J}_{\phi_0, D_0}$
 is the identity map for any ${\bf J}\in{\cal J}(P_{\phi_0},\Omega_{\phi_0})$.}

 Now if a loop $\chi_{t\in [0,1]}$ is homotopic to $\phi_0$ in ${\rm Ham}^c_0(V,\omega)$
 then there exists a  Hamiltonian fibre bundle
 isomorphism ${\rm I}^{\phi_0\chi}$ having {\rm CS} property from
 $P_{\phi_0}$ to $P_\chi$. We call $\Gamma_\omega$-
 equivalence class ${\rm I}^{\phi_0\chi}_\ast(D_0)$ of sections of $P_\chi$ as
 the {\it trivial class}. It is independent of choice of the isomorphism
 ${\rm I}^{\phi_0\chi}$ having {\rm CS} property.
 Thus $\Psi^{\bf J}_{\chi, T}$ is the identity map for
  the trivial class $T$ and any ${\bf J}\in{\cal J}(P_{\chi},\Omega_{\chi})$.

  As done in [LMP] the key point of the proof of Theorem 6.4 is to prove the
 composition rule for maps $\Psi_{\phi,D}$. This needs us to consider the relation
 between $P_\phi$, $P_\phi$ and $P_{\psi\ast\phi}$.
 However, unlike  the case of [LMP] under which there is the only coupling class
 $u_{\phi}$ corresponding to $\phi$, in our case we need to replace it by
 a suitable thing.
 For two smooth loops $\phi_{t\in [0,1]}$ and $\psi_{t\in [0,1]}$
 in ${\rm Ham}^c_0(V,\omega)$ we make the following assumptions:
 for a fixed sufficiently small $\epsilon>0$
$\phi_t=id$ for $t\notin [1/2+\epsilon, 1-\epsilon]$ and
$\psi_t=id$ for $t\notin [\epsilon, 1/2-\epsilon]$. Notice that they have been
extended to $\mbox{\Bb R}$ $1$-periodically. Let $H_\phi:S^1\times V\to\mbox{\Bb R}$
and  $H_\psi:S^1\times V\to\mbox{\Bb R}$ be the functions with compact support
and generating loops $\phi_{t\in [0,1]}$ and $\psi_{t\in [0,1]}$ respectively.
We can require them to satisfy: $H_\phi(t,\cdot)=0$ for
$t\notin [1/2+\epsilon, 1-\epsilon]$ and $H_\psi(t,\cdot)=0$
for $t\notin [\epsilon, 1/2-\epsilon]$.
Denote by
\begin{eqnarray*}
&&T_\epsilon:=\{[r,t]_+, [r,t]_-\in S^2\,|\,2/3\le r\le 1,\; t\in [1/2+\epsilon, 1-\epsilon]\},\\
&&T^\ast_\epsilon:=\{[r,t]_+, [r,t]_-\in S^2\,|\,2/3\le r\le 1,\; t\in [\epsilon, 1/2-\epsilon]\},\\
&&S^2_+:=\{[r,t]_+, [r,t]_-\in S^2\,|\, 0\le r\le 1,\; 0< t< 1/2\},\\
&&S^2_-:=\{[r,t]_+, [r,t]_-\in S^2\,|\, 0\le r\le 1,\; 1/2< t< 1\}.
\end{eqnarray*}
Clearly, $T_\epsilon$ and $T^\ast_\epsilon$ are proper subsets of the open left hemisphere
$S^2_+$ and open right hemisphere $S^2_-$ respectively.
>From the previous construction we may know that
$$P_\phi|_{S^2\setminus T_\epsilon}= (S^2\setminus T_\epsilon)\times V\qquad{\rm and}\qquad
P_\psi|_{S^2\setminus T^\ast_\epsilon}= (S^2\setminus
T^\ast_\epsilon)\times V.\leqno(59)$$ Thus we may construct the
fibre  sum $P_\phi\sharp P_\psi$ as follows: gluing
$P_\phi|_{S^2\setminus S^2_+}$ and $P_\psi|_{S^2\setminus S^2_-}$
along
$$\partial P_\phi|_{S^2\setminus S^2_+}=\{[r,t]_+, [r,t]_-\in S^2\,|\, t=0, 1/2\}\times V=
\partial P_\phi|_{S^2\setminus S^2_-}\leqno(60)$$
by the map: $[r, 0, x]^\phi_\pm\to [r, 0, x]^\psi_\pm$, and
$[r, 1/2, x]^\phi_\pm\to [r, 1/2, x]^\psi_\pm$.

On the other hand it is easy to know that under our assumptions
the composite loop $(\phi\ast\psi)_{t\in [0,1]}=(\phi_t\circ\psi_t)_{t\in [0,1]}$
is generated by the Hamiltonian function
$H_{\phi\ast\psi}:S^1\times V\to\mbox{\Bb R}$ given by
$$H_{\phi\ast\psi}(t,x)=\cases{H_\psi(t,x), &if $0\le t\le 1/2$;\cr
H_\phi(t,x), &if $1/2\le t\le 1$.\cr}\leqno(61)$$ Notice that the
Hamiltonian forms $\widetilde\Omega_\phi$ on $P_\phi$ and
$\widetilde\Omega_\psi$ on $P_\psi$ constructed as before
satisfies:
$$\widetilde\Omega_\phi|P_\phi|_{S^2\setminus T_\epsilon}=p_2^\ast\omega,\qquad
 \widetilde\Omega_\psi|P_\psi|_{S^2\setminus T^\ast_\epsilon}=p_2^\ast\omega.\leqno(62)$$
Hence under the fibre  sum operation they define a closed two-form $P_\phi\sharp P_\psi$,
denoted by $\widetilde\Omega_\phi\sharp\widetilde\Omega_\psi$. From the above construction
it is easily checked that $P_\phi\sharp P_\psi=P_{\phi\ast\psi}$ and the
closed two-form $\widetilde\Omega_\phi\sharp\widetilde\Omega_\psi$
 is exactly $\widetilde\Omega_{\phi\ast\psi}$ constructed
in the previous way, that is,
$$\widetilde\Omega_\phi\sharp\widetilde\Omega_\psi=\widetilde\Omega_{\phi\ast\psi}.\leqno(63)$$
Now for given sections $s$ of $P_\phi$ and $s^\prime$ of $P_\psi$,
by the section homotopy we assume that the restriction of $s$ to
$S^2\setminus T_{\epsilon/2}$ and that of $s^\prime$ to
$S^2\setminus T^\ast_{\epsilon/2}$ have the following versions
respectively,
$$s(z)=(z, v_0),\, z\in S^2\setminus T_{\epsilon/2},\quad{\rm and}\quad
s^\prime(z)=(z, v_0),\, z\in S^2\setminus
T^\ast_{\epsilon/2}\leqno(64)$$ for some fixed $v_0\in V$. Hence
they fit together to give one section of the bundle $P_\phi\sharp
P_\psi$, denoted by
$$s\sharp s^\prime.\leqno(65)$$
Combing (63) with (65) we get
$$\widetilde\Omega_\phi(s) + \widetilde\Omega_\psi(s^\prime)=
\widetilde\Omega_\phi\sharp\widetilde\Omega_\psi(s\sharp
s^\prime)= \widetilde\Omega_{\phi\ast\psi}(s\sharp
s^\prime).\leqno(66)$$ For such chosen sections $s$ and $s^\prime$
it follows from (59) that
$$c_1(TP_\phi^{\rm vert})(s)+ c_1(TP_\psi^{\rm vert})(s^\prime)=
c_1(TP_{\phi\ast\psi}^{\rm vert})(s\sharp s^\prime).\leqno(67)$$
In fact, since $c_1(TP_\phi^{\rm vert})(s)=c_1(s^\ast TP_\phi^{\rm
vert})([S^2])$, by the well-known Splitting Principle we only need
to consider the case of complex line bundle on $S^2$. The latter
case may be directly proved with Theorem 2.71 in [McSa2].

Notice that (66) and (67) lead to a natural map from
$$\Gamma_\omega(P_\phi)\times\Gamma_\omega(P_\psi)\to\Gamma_\omega(P_{\phi\ast\psi}),\;
(D, D^\prime)\to D\sharp D^\prime,\leqno(68)$$ where
$\Gamma_\omega(P_\phi)$, $\Gamma_\omega(P_\psi)$ and
$\Gamma_\omega(P_{\phi\ast\psi})$ are the sets of
$\Gamma_\omega$-equivalence classes of the sections of the bundles
$P_\phi$, $P_\psi$ and $P_{\phi\ast\psi}$ respectively. Similarly,
since ${\cal J}(P_\phi,\Omega_\phi)$ and ${\cal J}(P_\psi,
\Omega_\psi)$ are contractible using (59) we always choose ${\bf
J}\in{\cal J}(P_\phi,\Omega_\phi)$ and ${\bf J}^\prime\in{\cal
J}(P_\psi,\Omega_\psi)$ such that
$$J_z=J=J^\prime_w,\;\forall z\in S^2\setminus T_{\epsilon/2}\;{\rm and}\;
\forall w\in S^2\setminus T^\ast_{\epsilon/2},\leqno(69)$$ where
$J\in{\cal J}(V,\omega)$ such that $(V,\omega, J, g)$ is g.bounded
for $g\in{\cal G}(V)$. Then ${\bf J}$ and ${\bf J}^\prime$ fit
together to give one element in ${\cal
J}(P_{\phi\ast\psi},\widetilde\Omega_{\phi\ast\psi})$, denoted by
${\bf J}\sharp {\bf J}^\prime$. What we wish to prove is the
following composition rule.

\noindent{\bf Proposition 6.7.}\hspace{2mm}{\it For any
$\Gamma_\omega$-equivalence classes $D$ of sections of $P_\phi$
and $D^\prime$ of sections of $P_\psi$ it holds that
$$\Psi^{{\bf J}^\prime}_{\psi, D^\prime}\circ\Psi^{\bf J}_{\phi, D}=
\Psi^{{\bf J}\sharp {\bf J}^\prime}_{\phi\ast\psi, D\sharp
D^\prime}.\leqno(70)$$}

Before giving its proof we make an notation:

\noindent{\bf Remark 6.8}\hspace{2mm}The above base spaces of $P_\phi$
and $P_\psi$ are denoted by $S^2_L$ and $S^2_R$ respectively. Moreover, when
constructing the fibre sum $P_\phi\sharp P_\psi$ we will glue
$P_\phi|_{S^2\setminus S^2_{L\varepsilon}}$ and
$P_\psi|_{S^2\setminus S^2_{R\varepsilon}}$ along boundaries
$\partial P_\phi|_{S^2\setminus S^2_{L\varepsilon}}$ and
$\partial P_\psi|_{S^2\setminus S^2_{R\varepsilon}}$ by the map
$$\Bigl[\frac{\cos\varepsilon}{\cos(t+\varepsilon)}, t, x\Bigr]^\phi\mapsto
\Bigl[\frac{\cos\varepsilon}{\cos (t-\pi)}, t-\pi,
x\Bigr]^\psi,\leqno(71)$$ where
\begin{eqnarray*}
&&S^2_{L\varepsilon}:=\Bigl\{[r,t]_+, [r,t]_-\in S^2\,\Bigm|\,
\frac{\cos\varepsilon}{\cos t}< r\le 1,\; -\varepsilon< t< \varepsilon\Bigr\},\\
&&S^2_{R\varepsilon}:=\Bigl\{[r,t]_+, [r,t]_-\in S^2\,\Bigm|\,
\frac{\cos\varepsilon}{\cos(t+\varepsilon)}< r\le 1,\; \pi-\varepsilon< t<\pi +\varepsilon\Bigr\}.
\end{eqnarray*}
We denote the fibre sum by $P_\phi\sharp_\varepsilon P_\psi$. Notice that there exists
the canonical fibre bundle isomorphism ${\bf I }_\varepsilon$ from
$P_\phi\sharp_\varepsilon P_\psi$ to  $P_\phi\sharp P_\psi$.
Later, when  saying $\widetilde\Omega_{\phi\ast\psi}$ on $P_\phi\sharp_\varepsilon P_\psi$
and $c_1(TP^{\rm vert}_{\phi\ast\psi})$ we always mean them to be
the pullback of $\widetilde\Omega_{\phi\ast\psi}$ and
$c_1(TP^{\rm vert}_{\phi\ast\psi})$ on $P_{\phi\ast\psi}$ under ${\bf I}_\varepsilon^\ast$
without special statements. The sum $s\sharp s^\prime$ of sections and other related objects
will be understood similarly.\hfill$\Box$

Denote by $d=c_1(TP^{\rm vert}_\phi)(D)$ and $d^\prime=c_1(TP^{\rm vert}_\phi)(D^\prime)$.
By (67) it holds that
$$c_1(TP_{\phi\ast\psi}^{\rm vert})(D\sharp D^\prime)=d+ d^\prime.\leqno(72)$$
Thus both
$\Psi^{{\bf J}^\prime}_{\psi, D^\prime}\circ\Psi^{\bf J}_{\phi, D}$ and
$\Psi^{{\bf J}\sharp {\bf J}^\prime}_{\phi\ast\psi, D\sharp D^\prime}$
are the homomorphisms from $QH_\ast(M)$ to $QH_{\ast+ d+ d^\prime}(M)$.
For a given $\alpha\in H_\ast(M,\mbox{\Bb Z})$ the straightforward computations
shows
$$\Psi^{{\bf J}^\prime}_{\psi, D^\prime}\circ\Psi^{\bf J}_{\phi, D}(\alpha)
=\sum_{A\in\Gamma_\omega}(\sum_{B\in\Gamma_\omega}\alpha_{B,
A-B})\otimes e^A,\leqno(73)$$ where $\alpha_{B,A-B}\in H_\ast(V)$
is determined by
$$\alpha_{B, A-B}\cdot_V\beta=\Phi_{(\psi,D^\prime+A-B)}(\alpha_B, \beta),\;\forall \beta\in H_\ast(V),\leqno(74)$$
and $\alpha_B\in H_\ast(V)$ by
$$\alpha_B\cdot_V\gamma=\Phi_{(\phi, D+B)}(\alpha,\gamma),\;\forall\gamma\in H_\ast(V).\leqno(75)$$
Notice that we also have
$$\dim\alpha_{B,A-B}=\dim\alpha_B+ 2c_1(TP_\psi^{\rm vert})(D^\prime)+2c_1(A-B),\leqno(76)$$\vspace{-3mm}
$$\dim\alpha_B=\dim\alpha+ 2c_1(TP_\phi^{\rm vert})(D)+2c_1(B).\leqno(77)$$
Moreover, by definition we also have
$$\Psi^{{\bf J}\sharp {\bf J}^\prime}_{\phi\ast\psi, D\sharp D^\prime}(\alpha)=\sum_{A\in\Gamma_\omega}
\hat\alpha_A\otimes e^A,\leqno(78)$$ where $\hat\alpha_A\in
H_\ast(V)$ is determined by
$$\hat\alpha_A\cdot_V\gamma=\Phi_{(\phi\ast\psi, D\sharp D^\prime+ A)}(\gamma),\;\gamma\in H_\ast(V),\leqno(79)$$
\vspace{-3mm}
$$\dim\hat\alpha_A=\dim\alpha+2c_1(TP_{\phi\ast\psi}^{\rm vert})(D\sharp D^\prime)+2c_1(A).\leqno(80)$$
Thus we only need to prove that
$$\hat\alpha_A=\sum_{B\in\Gamma_\omega}\alpha_{B,A-B},\;\forall A\in \Gamma_\omega.\leqno(81)$$

To complete the proof of Proposition 6.7 we need several lemmas.

\noindent{\bf Lemma 6.9.}\hspace{2mm}{\it For every fixed $A\in
H_\ast(V)$ the sum of right side in (81) is always finite.}

\noindent{\bf Lemma 6.10.}\hspace{2mm}{\it There exist the regular
almost complex structures $\hat J$ on $P_\phi$ and $\hat J^\prime$
on $P_\psi$ such that they agree on gluing domain of
$P_{\phi}\sharp_\varepsilon P_\psi$.}

Without special statements we shall fix $\hat J$ and $\hat J^\prime$.
The proof of Lemma 6.9 is given after Lemma 6.12 and  Lemma 6.10
will be proved at the end of this section.

Following the notations in \S4.

\noindent{\bf Lemma 6.11.}\hspace{2mm}{\it Let $e_1:U\to V$ and
$e_2:U\to V$ be two $C^p$-smooth pseudo-cycles, and $\alpha:A\to
V$ and $\beta:B\to V$ be two $C^q$-smooth pseudo-cycles. Assume
that
$$\dim U+\dim A + \dim B\ge 2\dim V,\leqno(83)$$
then for every sufficiently large integer $r>{\rm min}\{p,q\}$
there exists a set
$\chi^r(V,e_1, e_2,\alpha,\beta)\subset\chi^r(V)_0\times\chi^r(V)_0$ of the second category
such that $e=(e_1, e_2)$ is transverse to
$({\cal F}^r({\rm X})\circ \alpha)\times ({\cal F}^r({\rm Y})\circ \beta)$ for all
$({\rm X}, {\rm Y})\in\chi^r(V,e_1, e_2,\alpha,\beta)$.
These $\chi^r(V,e_1, e_2, \alpha,\beta)$ also satisfy:
$$\chi^r(V,e_1, e_2,\alpha,\beta)\supseteq\chi^{r+1}(V,e_1, e_2,\alpha,\beta)\supseteq\cdots,$$
which implies that for any $({\rm X},{\rm Y})\in\chi^r(V,e_1,e_2,\alpha,\beta)$ and
$({\rm X}^\prime,{\rm Y}^\prime)\in\chi^s(V,e_1,e_2,\alpha,\beta)$ with $s>r$ it holds that
$$[({\cal F}^r({\rm X})\circ \alpha)\times ({\cal F}^r({\rm Y})\circ \beta)]\cdot e=
[({\cal F}^s({\rm X}^\prime)\circ\alpha\times ({\cal F}^s({\rm Y}^\prime)]\cdot e$$
provided that the equality in (83) holds, and one of $e=(e_1, e_2)$ and $\alpha\times\beta$
is the strong pseudo-cycle.}

The proof of this lemma is similar to that of Lemma 4.1.
Replacing (27) one only consider the map
$$\Xi:U\times A\times B\times\chi^r(V)_0\times\chi^r(V)_0\to V\times V\times V\times V$$
given by
$$(u,a,b,{\rm X},{\rm Y})\mapsto\Bigl((e_1(u), e_2(u)),
\bigl({\cal F}^r({\rm X})(\alpha(a)), {\cal F}^r({\rm Y})(\beta(b))\bigr)\Bigr).$$
It is easy to prove that it is transverse to
$$\Delta_{V\times V}:=\{(u,v, u,v)\,|\,(u,v)\in V\times V\}.$$
The standard arguments may finish the proof.

By (i) of Lemma 4.1 one know that if
$$\dim U+\dim A\ge\dim V\leqno(84)$$
then  for every sufficiently large integer $r>{\rm min}\{p,q\}$
there exists a set
$\chi^r(V,e_1,\alpha)\subset\chi^r(V)_0$ of the second category
such that $e_1$ is transverse to
${\cal F}^r({\rm X})\circ \alpha$ for all
${\rm X}\in\chi^r(V,e_1,\alpha)$. From Claim A.1.11 of [LeO] the space
$$\chi^r_1(V,e_1\times e_2,\alpha,\beta)\leqno(85)$$
consisting of all ${\rm X}\in\chi^r(V)_0$ for which the intersection
$$\chi^r(V,e_1, e_2,\alpha,\beta)\cap[\{\rm X\}\times\chi^r(V)_0]$$
is a countable intersection of open dense subset in
$\{\rm X\}\times\chi^r(V)_0$ must be a countable
intersection of open dense subsets in $\chi^r(V)_0$.
Thus the intersection
$$\chi^r_1(V,e_1\times e_2,\alpha,\beta)\cap\chi^r(V,e_1,\alpha)\leqno(86)$$
is also a countable
intersection of open dense subsets in $\chi^r(V)_0$.
For every ${\rm X}$ in this intersection there must be a ${\rm Y}\in\chi^r(V)_0$
such that $({\rm X},{\rm Y})\in\chi^r(V,e_1\times e_2,\alpha,\beta)$. Thus this
pair $({\rm X},{\rm Y})$ satisfies:
\begin{description}
\item[(i)]$e=(e_1, e_2)$ is transverse to
$({\cal F}^r({\rm X})\circ \alpha)\times ({\cal F}^r({\rm Y})\circ \beta)$,
\item[(ii)]$e_1$ is transverse to ${\cal F}^r({\rm X})\circ \alpha$
\end{description}
under the assumptions (83) (84).

\noindent{\bf Lemma 6.12.}\hspace{2mm}{\it Let $e_i:U\to V, i=1,2$
and $\alpha:A\to V$ and $\beta:B\to V$ be all $C^r$-smooth
pseudo-cycles. Assume that $e=(e_1,e_2)$ is transverse to
$\alpha\times\beta$, $e_1$ is transverse to $\alpha$ and (83) (84)
also hold. Then
$$\Delta(U\times A):=\{(u,a)\,|\,e_1(u)=\alpha(a)\}\leqno(87)$$
is a $C^r$-smooth manifold of dimension
$\dim U+\dim A-\dim V$, and
$$\hat e_2:\Delta(U\times A)\to V,\;(u,a)\mapsto e_2(u)\leqno(88)$$
is also $C^r$-smooth pseudo-cycle which is transverse to
$\beta$. Moreover, if $\alpha$ and $\beta$ are strong
pseudo-cycle then it holds that
$$e\cdot(\alpha\times\beta)=\hat e_2\cdot\beta.\leqno(89)$$}

\noindent{\it Proof}.\hspace{2mm}Let $\hat e_2(u,a)=\beta(b)$. We
wish to prove
$$D\hat e_2(u,a)(T_{(u,a)}\Delta(U\times A))+ D\beta(b)(T_b B)=
T_{\beta(b)}V.$$
Notice that
$$T_{(u,a)}\Delta(U\times A)=\{(\vec u,\vec a)\in T_u U\times T_a A\,|\,
De_1(u)(\vec u)=D\alpha(a)(\vec a)\}.$$
It suffice to prove that for any $\xi\in T_{\beta(b)}V$ there
exist $\vec u\in T_u U$, $\vec a\in T_a A$ and $\vec b\in T_b B$
such that
$$De_1(u)(\vec u)=D\alpha(a)(\vec a)\quad{\rm and}\quad
De_2(u)(\vec u)+D\beta(b)(\vec b)=\xi.\leqno(90)$$ But $e$ is
transverse to $\alpha\times\beta$. Therefore, there exist $(\vec
u,\vec a,\vec b)\in T_u U\times T_a A\times T_b V$ such that
$$De(u)(\vec u)+D(\alpha\times\beta)(a,b)(-\vec a,\vec b)=(0,\xi).$$
Clearly, this is equivalent to (90). By similar arguments for the
boundary parts we can prove that $\hat e_2$ is a $C^r$-smooth
pseudo-cycle which is transverse to $\beta$.

Notice that $\hat e_2$ is also a strong pseudo-cycle if $\alpha$
is. Now $e\cdot(\alpha\times\beta)$ and $\hat e_2\cdot\beta$ are
well-defined. To prove them being equal we notice that
\begin{eqnarray*}
\hat e_2\cdot\beta&=&\sum_{\scriptstyle e_1(u)=\alpha(a)\atop\scriptstyle
e_2(u)=\beta(b)}{\it sign}\Bigl((u,a),b\Bigr)\\
&=&\sum_{\scriptstyle e_1(u)=\alpha(a)\atop\scriptstyle
e_2(u)=\beta(b)}{\it sign}(u,a,b),\\
e\cdot(\alpha\times\beta)&=&\sum_{\scriptstyle e(u)=(\alpha\times\beta)(a,b)}
{\it sign}\Bigl(u, (a,b)\Bigr)\\
&=&\sum_{\scriptstyle e_1(u)=\alpha(a)\atop\scriptstyle
e_2(u)=\beta(b)}{\it sign}(u,a,b).
\end{eqnarray*}
Here some details on the orientation are omitted. It is not very difficult
to give them. At least, for the mod 2 intersection number the above arguments
is completed.
The lemma is proved.\hfill$\Box$\vspace{2mm}

\noindent{\bf Remark 6.13.}\hspace{2mm} Using Lemma 6.12 we may
give a pseudo-cycle expression of $\alpha_B$ in (75) as follows:
Firstly, by Lemma 6.11 $H\in{\rm Diff}(V\times V)$ in (52) may be
chosen as the form $H=(h_1, h_2)$ with $h_i\in{\rm Diff}(V)$,
$i=1,2$. Thus (75) becomes
$$\alpha_B\cdot_V\gamma={\rm EV}^{\phi(D+B)}\cdot(h_1\circ f_M\times h_2\circ f_N),\leqno(91)$$
where $f_M:M\to V$ and $f_N:N\to V$ are the strong pseudo-cycle representatives of $\alpha$
and $\gamma$ respectively, $h_i\in{\rm Diff}(V)$, $i=1,2$, and
$${\rm EV}^{\phi(D+B)}=({\rm EV}_1^{\phi(D+B)}, {\rm EV}_2^{\phi(D+B)}):
{\cal S}(P_\phi,\Omega_\phi,j,\hat J, D+B)\to V\times V,\;
s\mapsto\bigr(F_1^\phi(s(z_1^\phi)), F_2^\phi(s(z_2^\phi))\bigl)$$
is the pseudo-cycle determined by the evaluation map. By lemma
6.12 the right side of (91) is equal to
$$\widehat{\rm EV}_2^{\phi(D+B)}\cdot (h_2\circ f_N),\leqno(92)$$
where the pseudo-cycle
$$\widehat{\rm EV}_2^{\phi(D+B)}:\Delta({\cal S}(P_\phi,\Omega_\phi,j,\hat J, D+B), M)\to V\leqno(93)$$
given by
$$(s, a)\mapsto F_2^\phi(s(z_2^\phi)).$$
By definition
$$\Delta({\cal S}(P_\phi,\Omega_\phi,j,\hat J, D+B), M)=\Bigl\{(s, a)\in
{\cal S}(P_\phi,\Omega_\phi,j,\hat J, D+B)\times M\,\Bigm|\, F_1^\phi(s(z_1^\phi))=h_1\circ f_M(a)\Bigl\}.$$
Thus (93) may be considered as a pseudo-cycle representative of $\alpha_B$.
\hfill$\Box$\vspace{2mm}

\noindent{\it Proof of Lemma 6.9}.\hspace{2mm}Assume that there
exists $A\in\Gamma_\omega$ such that $\alpha_{B, A-B}\ne 0$ for
infinitely many $B\in\Gamma_\omega$. Denote them by $B_1,
B_2,\cdots$. By Remark 6.13 one gets infinitely many pseudo-cycles
$$\widehat{\rm EV}_2^{\phi(D+B_i)}:\Delta({\cal S}(P_\phi,\Omega_\phi,j,\hat J, D+B_i), M)\to V,\leqno(94)$$
\vspace{-4mm}
$$\Delta({\cal S}(P_\phi,\Omega_\phi,j,\hat J, D+B_i), M)=\Bigl\{(s, a)\in
{\cal S}(P_\phi,\Omega_\phi,j,\hat J, D+B_i)\times M\,|\, F_1^\phi(s(z_1^\phi))=h^{(i)}_1\circ f_M(a)\Bigl\}$$
for some $h^{(i)}_1\in{\rm Diff}(V)$.
>From Lemma 6.11 and the arguments under it one can assume all $h_1^{(i)}$ to be the same $h_1$.
But the image of $h_1\circ f_M$ is contained in a compact subset of $V$. From the results in \S2 it follows
that the image sets of all sections $s$ which are such that
$$\Bigl(\{s\}\times M\Bigr)\bigcap\biggl(\bigcup\Delta({\cal S}(P_\phi,\Omega_\phi,j,\hat J, D+B_i), M)\biggr)\ne
\emptyset,\leqno(95)$$ are contained in a compact subset of
$P_\phi$. Thus the image sets of all such pseudo-cycles
representatives of $\alpha_{B_i}$ given by Remark 6.13 are
contained in a compact subset $K(\phi)$ of $V$. By the assumption
at the beginning
$$\alpha_{B_i, A-B_i}\cdot_V\beta_i\ne 0,\;{\it for\;some}\;\beta_i\in H_\ast(V)\leqno(96)$$
Now from (74) it follows that there exist sections
$$s_i^\prime\in{\cal S}(P_\psi,\Omega_\psi,j,\hat J^\prime, D^\prime+A-B_i)$$
such that ${\rm EV}_1^{\psi(D^\prime+A-B_i)}(s_i^\prime)=F^\psi_1(s^\prime_i(z_1^\psi))$
are contained in the compact subset $K(\phi)$.
Hence the image sets of all sections $s^\prime_i$ are contained in a compact subset
${\rm S}(\psi)$ of $P_\psi$. Because all $s^\prime_i$ are $(j,\hat J^\prime)$-holomorphic
it holds that
$$(\widetilde\Omega_\psi+ c_0\sigma)(s^\prime_i)\ge 0,\;i=1,2,\cdots,$$
which implies
$$\omega(B_i)=(\widetilde\Omega_\psi+ c_0\sigma)(B_i)\le
(\widetilde\Omega_\psi+ c_0\sigma)(D^\prime+A),\;i=1,2,\cdots.
\leqno(97)$$ Hence
$$(\widetilde\Omega_\phi+c_0\sigma)(D+B_i)\le (\widetilde\Omega_\phi+c_0\sigma)(D)+
(\widetilde\Omega_\psi+
c_0\sigma)(D^\prime+A),\;i=1,2,\cdots.\leqno(98)$$ Take $(s_i,
a_i)\in\Delta({\cal S}(P_\phi,\Omega_\phi,j,\hat J, D+B_i), M)$
one gets infinitely many $(j,\hat J)$-holomorphic sections
$\{s_i\}$ which represent infinitely many different classes and
 whose image sets are contained
in a fixed compact subset of a g.bounded symplectic manifold
$(P_\phi,\widetilde\Omega_\phi+ c_0\sigma, \hat J, g)$. With the same reason as in Lemma 6.5
(98) leads to a contradiction.
\hfill$\Box$\vspace{2mm}

Now we have known that the sum on the right side of (81) is
actually finite sum. To prove (81) holding let us check their
pseudo-cycle representatives given by Remark 6.13. Using the
pseudo-cycle representative of $\alpha_B$ given by (93) one may
get a pseudo-cycle representative of $\alpha_{B, A-B}$ as follows:
$$\widehat{\rm EV}_2^{\psi(D^\prime+A-B)}:
\Delta\biggl({\cal S}(P_\psi,\Omega_\psi, j,\hat J^\prime,
D^\prime+A-B), \Delta({\cal S}(P_\phi,\Omega_\phi, j,\hat J, D+B),
M)\biggr)\to V\leqno(99)$$ given by
$$\Bigl(s^\prime, (s,a)\Bigr)\mapsto F^\psi_2(s^\prime(z^\psi_2)).\leqno(100)$$
By definition it is easy to check that the set in (99) consists of all
triples $(s^\prime, s, a)$ satisfying the conditions
$$\left.\matrix{&s^\prime\in {\cal S}(P_\psi,\Omega_\psi, j,\hat J^\prime, D^\prime+A-B)\cr
&s\in{\cal S}(P_\phi,\Omega_\phi, j,\hat J, D+B)\cr &a\in M\cr
&F_2^\phi (s(z_2^\phi))=h_1^{\phi (D+B)}\circ f_M(a)\cr &F_1^\psi
(s^\prime(z_1^\psi))=h_1^{\psi (D^\prime+A-B)}\circ F_2^\phi
(s(z_2^\phi))} \right\}\leqno(101)$$ for some $h_1^{\phi(D+B)}$
and $h_1^{\psi(D^\prime+A-B)}$ in ${\rm Diff}(V)$. Moreover, from
Lemma 6.12 it is easily computed that the dimension of manifold in
(99) is
$$\dim\alpha+ 2c_1(A) + 2c_1(TP_{\phi\ast\psi}^{\rm vert})(D\sharp D^\prime).\leqno(102)$$

On the other hand $\hat\alpha_A$ has the pseudo-cycle representative:
$$\widehat{\rm EV}_2^{\phi\ast\psi(D\sharp D^\prime+A)}:
\Delta\Bigl({\cal S}(P_{\phi\ast\psi},\Omega_{\phi\ast\psi}, j,
\hat J\sharp\hat J^\prime, D\sharp D^\prime +A), M\Bigr)\to
V\leqno(103)$$ given by
$$(\sigma, a)\mapsto F_2^{\phi\ast\psi}(\sigma(z_2^{\phi\ast\psi})).\leqno(104)$$
By definition
$$\Delta\Bigl({\cal S}(P_{\phi\ast\psi},\Omega_{\phi\ast\psi}, j, \hat J\sharp\hat J^\prime,
D\sharp D^\prime +A), M\Bigr)\leqno(105)$$ consists of all pairs
$(\sigma,a)$ satisfying
$$\left.\matrix{&\sigma\in{\cal S}(P_{\phi\ast\psi},\Omega_{\phi\ast\psi}, j, \hat J\sharp\hat J^\prime,
 D\sharp D^\prime +A)\cr
&a\in M\cr
&F_1^{\phi\ast\psi}(\sigma(z_1^{\phi\ast\psi}))=h_1^{\phi\ast\psi(D\sharp
D^\prime+A)}\circ f_M(a) }\right\}\leqno(106)$$ for some
$h_1^{\phi\ast\psi(D\sharp D^\prime+A)}$ in ${\rm Diff}(M)$. Here
it should be noted that the choices of $h_1^{\phi\ast\psi(D\sharp
D^\prime+A)}$ in (106) and $h_1^{\phi (D+B)}$ and $h_1^{\psi
(D^\prime+A-B)}$ in (101) have the ``bigger'' freedom. But the
choices of $h_1^{\psi (D^\prime+A-B)}$ are under the case that
$h_1^{\phi (D+B)}$ is chosen. Another important point is the maps
in (99) and (103) to have precompact image sets in $V$. Thus they
are all strong pseudo-cycles in the sense of \S4.

Having the above preparation we may prove (81) and thus finish the proof of Proposition 6.7.
We only need to prove
$$PD(\hat\alpha_A)=\sum_{B\in\Gamma_\omega}PD(\alpha_{B,A-B}),\;\forall A\in\Gamma_\omega.\leqno(107)$$
That is, their Poincar\`e dualities in $H^\ast_c(V)$ are same. But
(107) is equivalent to
$$\langle PD(\hat\alpha_A), \gamma\rangle=\sum_{B\in\Gamma_\omega}
\langle PD(\alpha_{B,A-B}),\gamma\rangle,\;\forall \gamma\in
H_\ast(V).\leqno(108)$$ Therefore, one only need to prove that for
every $\gamma\in H_\ast(V)$ with
$$\dim\gamma=\dim\alpha+ 2c_1(A) + 2c_1(TP_{\phi\ast\psi}^{\rm vert})(D\sharp D^\prime)$$
we may choose a pseudo-cycle representative of it $\Upsilon:T\to V$ such that
it is transverse to the map in (104) and all maps in (99) and
$$\Upsilon\cdot\widehat{\rm EV}_2^{\phi\ast\psi(D\sharp D^\prime+A)}=
\sum_{B\in\Gamma_\omega}\Upsilon\cdot\widehat{\rm
EV}_2^{\psi(D^\prime+A-B)}.\leqno(109)$$ By definitions the left
side of (109) is equal to the sum
$$\sum{\rm sign}(r,\sigma, a)\leqno(110)$$
when $(r,\sigma,a)$ takes over the set
$$\Biggl\{\bigl(r, (\sigma, a)\bigr)\in T\times\Delta\Bigl({\cal S}(P_{\phi\ast\psi},
\Omega_{\phi\ast\psi}, j, \hat J\sharp\hat J^\prime, D\sharp
D^\prime +A), M\Bigr)
\,\Biggm|\,\Upsilon(r)=F_2^{\phi\ast\psi}(\sigma(z_2^{\phi\ast\psi}))\Biggr\}.\leqno(111)$$
The right side of (109) is equal to the sum
$$\sum{\rm sign}(r, s^\prime, s, a),\leqno(112)$$
where $\biggl(r,\Bigl(s^\prime,(s,a)\Bigr)\biggr)$ takes over the set
$$\Lambda(T, M, P_\psi, P_\phi, \hat J, \hat J^\prime, A, D, D^\prime)\leqno(113)$$
consisting of all
$$\biggl(r,\Bigl(s^\prime,(s,a)\Bigr)\biggr)\in T\times\bigcup_{B\in\Gamma_\omega}
\Delta\biggl({\cal S}(P_\psi,\Omega_\psi, j,\hat J^\prime,
D^\prime+A-B), \Delta({\cal S}(P_\phi,\Omega_\phi, j,\hat J, D+B),
M)\biggr)\leqno(114)$$ such that
$\Upsilon(r)=F^\psi_2(s^\prime(z^\psi_2))$. Notice that two sets
in (111) and (113) are finite.

By Remark 6.8 we here may choose
$$\left.\matrix{&z_1^\phi=[\frac{5}{6},\frac{1}{2}]_+=[\frac{5}{6},-\frac{1}{2}]_-,\quad
z_2^\phi=[\frac{5}{6},0]_+=[\frac{5}{6}, 0]_-\cr
&z_1^\psi=[\frac{5}{6},\frac{1}{2}]_+=[\frac{5}{6},-\frac{1}{2}]_-,\quad
z_2^\psi=[\frac{5}{6},0]_+=[\frac{5}{6},
0]_-}\right\}.\leqno(115)$$
Since the bundle $P_\psi$ and $P_\psi$
are trivial near $z_2^\phi$ and $z_1^\psi$ respectively, one can
use the gluing techniques developed in [RT1][McSa1] to prove that
there exists an orientation-preserving bijection between the set
in (111) and one in (113). This can lead to (109). Hence the proof
of Proposition 6.7 is completed under the assumption that Lemma
6.10 holds. \hfill$\Box$\vspace{2mm}

\noindent{\it Proof of Lemma 6.10}.\hspace{2mm}Recall the
technique used in \S2 and \S4. We only need to prove the following
fact:

\noindent{\bf Fact 6.14.}\hspace{2mm}{\it For a Riemannian vector
bundle $\pi:E\to W$, denoted $C^0_b(E)$ by the Banach space of all
bounded continuous sections of $\pi$. A norm of a section $s\in
C^0_b(E)$ is given by
$$\|s\|:=\sup_{x\in W}\|s(x)\|_g,$$
where $g$ is a given Riemannian metric on $E$.
Let $W_0$ be an open submanifold of $W$. Then for every open dense subset
${\cal A}$ in $C^0_b(E)$ the restriction ${\cal A}|W_0:=\{s|_{W_0}\,|\,s\in{\cal A}\}$
is also an open dense subset in $C^0_b(E|W_0)$.}

In fact, if there exists an open ball $B(s_0, \delta)\subseteq C^0_b(E|W_0)\setminus{\cal A}|W_0$
then one can find a section $s\in C^0_b(E)$ such that
$$\|s|_{W_0}-s_0\|<\frac{1}{5}\delta.$$
For this section $s$ there exists a section $s^\prime\in{\cal A}$ such that
$$\|s-s^\prime\|<\frac{1}{5}\delta.$$
Specially, this shows that $\|s|_{W_0}-s^\prime|_{W_0}\|<\frac{1}{5}\delta$.
Thus $s^\prime|_{W_0}\notin {\cal A}|W_0$, which leads to a contradiction.
\hfill$\Box$\vspace{2mm}

Now as in [LMP] it follows from Lemma 6.6 and Proposition 6.7 that
every $\Psi^{\bf J}_{\phi, D}$ is an isomorphism which leads to Theorem 6.4.

\noindent{\bf Remark 6.15.}\hspace{2mm}The conclusion of Corollary
6.3 can be actually strengthened to general case, that is, ${\rm
Diff}(M,\partial M)$ is replaced by ${\rm Diff}(M)$. We will
outline these as follows. Let $[0, 1]\to\phi_t$ be a smooth loop
in ${\rm Symp}(M,\omega)$, and $(\widetilde M,\widetilde\omega)$ a
noncompact symplectic manifold associated to $(M,\omega)$ as in
\S5. Here we need to write it in detail. Since $\partial M$ is a
hypersurface of contact type, for a contact form $\alpha$ on
$\partial M$ with $d\alpha=\omega|_{\partial M}$ the standard
arguments shows that there exists a $\varepsilon\in (0, 1)$ and an
embedding $\varphi:\partial M\times [\varepsilon, 1]\to M$ of
codimension zero such that
$$\varphi(m,1)=m\quad{\rm and}\quad\varphi^\ast\omega=d\Theta\;{\rm on}\;
\partial M\times [\varepsilon, 1],\leqno(116)$$
where $\Theta$ is a one-form on $\partial M\times [\varepsilon, +\infty)$
with $\Theta(m,z)=z\alpha(m)$ at a point $(m,z)$.
Then $(\widetilde M,\widetilde\omega)$ can be obtained by gluing
$(M,\omega)$ and $(\partial M\times [\varepsilon, +\infty), d\Theta)$
with $\varphi$. That is, $(m, z)\in \partial M\times [\varepsilon, 1]$
and $\varphi(m,z)\in M$ are identified. Notice that
$\phi_t(\partial M)=\partial M$, one can always find
a $\epsilon\in (\varepsilon, 1)$ such that
$$\bigcup_{t\in [0,1]}\phi_t\circ\varphi(\partial M\times (\epsilon, 1])$$
is contained in ${\rm Im}(\varphi)$. Thus every
$$\varphi^{-1}\circ\phi_t\circ\varphi:\partial M\times (\epsilon, 1]\to
 \partial M\times (\varepsilon, 1]$$
 is an embedding of codimension zero, and it also holds that
 $$\varphi^{-1}\circ\phi_t\circ\varphi(m, 1)=(\Phi_t(m), 1),\quad\forall m\in\partial M,$$
where $\Phi_t:\partial M\to\partial M$ is a smooth family of diffeomorphisms.
Since $(\varphi^{-1}\circ\phi_t\circ\varphi)^\ast d\Theta=d\Theta$ it must holds that
$$\Phi_t^\ast\alpha=\alpha.\leqno(117)$$
Define
$$\widetilde\phi_t: \widetilde M\to\widetilde M,\;q\mapsto
\cases{\phi_t(q) &on $q\in M$;\cr (\Phi_t(m), z)&on
$q=(m,z)\in\partial M\times [1,+\infty)$.\cr}\leqno(118)$$ It is
easily checked that $t\mapsto\widetilde\phi_t$ is a smooth loop in
${\rm Symp}(\widetilde M, \widetilde\omega)$. Moreover, if
$\{\phi_t\}_{t\in [0,1]}$ is generated by a smooth function
$H:M\times\mbox{\Bb R}/\mbox{\Bb Z}\to\mbox{\Bb R}$ then
$\{\widetilde\phi_t\}$ is generated by the smooth function
$$\widetilde H: \widetilde M\times\mbox{\Bb R}/\mbox{\Bb Z}\to\mbox{\Bb R},\;
(q, t)\mapsto\cases{ H(m, t) & if $(q,t)=(m, t)\in
M\times\mbox{\Bb R}/\mbox{\Bb Z}$;\cr H(m, t) &if $(q,t)=((m, z),
t)\in(\partial M\times [1,+\infty))\times\mbox{\Bb R}/\mbox{\Bb
Z}$.\cr}\leqno(119)$$ Now one may construct a Hamiltonian fibre
bundle $P_{\widetilde\phi}$ over $S^2$ with fibre $(\widetilde
M,\widetilde\omega)$. Furthermore, replacing $H$ with $\widetilde
H$ in the previous construction we may get a Hamiltonian $2$-form
$\widetilde\Omega_{\widetilde\phi}$ on $P_{\widetilde\phi}$. An
important point is that
$(P_{\widetilde\phi},\widetilde\Omega_{\widetilde\phi}, \hat J,
G)$ is also g.bounded for some $\hat J\in\widehat{\cal J}(j,{\bf
J})$ and some complete Riemannian metric $G$. Suitably modifying
the above arguments one may obtain the following corresponding
results to Theorem 6.4 and Theorem 6.1.

\noindent{\bf Proposition 6.16.}\hspace{2mm}{\it For a loop
$\phi_{t\in [0,1]}$ in ${\rm Ham}(M,\omega)$ and the extension
loop $\widetilde\phi_{t\in [0,1]}$ in ${\rm Ham}(\widetilde
M,\widetilde\omega)$ as above, the homomorphism
$i:H_\ast(\widetilde M,\mbox{\Bb Q})\to
H_\ast(P_{\widetilde\phi},\mbox{\Bb Q})$ is injective.
Consequently, the endomorphism $\partial_{\widetilde\phi}:
H_\ast(\widetilde M,\mbox{\Bb Q})\to H_{\ast+1}(\widetilde
M,\mbox{\Bb Q})$ vanishes. Especially, the endomorphism
$\partial_{\phi}: H_\ast(M,\mbox{\Bb Q})\to H_{\ast+1}(M,\mbox{\Bb
Q})$ vanishes.}

Using this result and the flux homomorphism theorem given in Appendix which is the version of Theorem 10.12 in
[McSa2] on the compact symplectic manifold with contact type boundary we get the following strengthened version
of Corollary 6.3.

\noindent{\bf Corollary 6.17.}\hspace{2mm}{\it For
$\phi\in\pi_1({\rm Diff}(M), id)$ and any two $\omega_1$ and
$\omega_2$ in ${\rm Cont}(M)$ it holds that $\phi\in {\rm
Im}(H_{\omega_1})\cap{\rm Im}(S_{\omega_2})$ if and only if
$\phi\in {\rm Im}(H_{\omega_2})\cap{\rm Im}(S_{\omega_1})$.}

Finally, we point out that using results in \S5 one can also generalize Theorem 5.A
in [LMP] to the present case.

\noindent{\bf Remark 6.18.}\hspace{2mm}After this paper had been
finished I saw D. McDuff's beautiful paper [Mc2]. It is very
possible to use our method to generalize her some results.
Moreover, from proof of Theorem 6.4 it easily follows that Theorem
6.4 still holds if the loop $\phi$ belongs to ${\rm
Ham}(V,\omega)$ rather than ${\rm Ham}^c(V,\omega)$, but we must
require that Hamiltonian function $H_\phi:S^1\times V\to\mbox{\Bb
R}$ generating $\phi$ satisfies some conditions( for example, a
possible choice is one that for some g.bounded Riemannian metric
$g$ on $V$ it holds that ${\rm Sup}\|dH(t,x)\|_g<+\infty$). These
will be given in other place.

\appendix
\section*{Appendix}

Suitably modifying the proof of Theorem 10.12 in [McSa2] one may get the following theorem.
For convenience of the readers we shall give its proof.

\noindent{\bf Theorem A.}\hspace{2mm}{\it Let $(M,\omega)$ be a
compact symplectic manifold with contact type boundary. Then a
smooth path
$$[0, 1]\to{\rm Symp}_0(M,\omega): t\mapsto\phi_t$$
from $\phi_0=id$ may be isotopic with fixed endpoints to a Hamiltonian path in ${\rm Ham}(M,\omega)$
if and only if ${\rm Flux}(\{\phi_t\})=0$.}

\noindent{\it Proof}.\hspace{2mm}Firstly, notice that  the flux
homomorphism is still
 well-defined on $\widetilde{\rm Symp}_0(M,\omega)$ or even on
$\widetilde{\rm Symp}_0(\widetilde M, \widetilde\omega)$ and
is indeed a homomorphism because there exists a natural homotopy equivalence between
$M$ and $\widetilde M$.

Next, we only need to prove the ``only if'' part.
Let $\phi_{t\in [0,1]}$ be a smooth path from $\phi_0=id$ in ${\rm Symp}_0(M,\omega)$ with
${\rm Flux}(\{\phi_t\})=0$. As in Remark 6.15 it is extended into a path from $id$
in ${\rm Symp}_0(\widetilde M,\widetilde\omega)$, denoted by $\widetilde\phi_{t\in [0,1]}$.
It has the version as in (117)(118).
Denote by
$$X_t=(\frac{d}{dt}\phi_t)\circ\phi_t^{-1},\quad
\widetilde
X_t=(\frac{d}{dt}\widetilde\phi_t)\circ\widetilde\phi_t^{-1}\quad{\rm
and}\quad {\cal
X}_t=(\frac{d}{dt}\Phi_t)\circ\Phi_t^{-1}\leqno(A.1)$$ then
$$\widetilde X_t(q)=\cases{ X_t(m) & if $q=m\in M$;\cr
({\cal X}_t(m), 0) &if $q=(m, z)\in\partial M\times
[1,+\infty)$.\cr}\leqno(A.2)$$ Thus
$$i_{\widetilde X_t}\widetilde\omega(q)=\cases{ i_{X_t}\omega(m) & if $q=m\in M$;\cr
-\alpha({\cal X}_t)(m)dz-zd(\alpha({\cal X}_t)(m) &if $q=(m,
z)\in\partial M\times [1,+\infty)$.\cr}\leqno(A.3)$$ Moreover, it
always holds that ${\rm Flux}(\{\phi_t\})={\rm
Flux}(\{\widetilde\phi_t\})$. Since ${\rm
Flux}(\{\widetilde\phi_t\})=0$, there exists a function
$\widetilde F:\widetilde M\to\mbox{\Bb R}$ such that
$$\int^1_0 i_{\widetilde X_t}\widetilde\omega dt=d\widetilde F.$$
It is easy to verify that up to a constant $\widetilde F|_{\partial M\times [1, +\infty)}$
may be chosen as:
$$\widetilde F(m,z)=-\int^1_0z\alpha({\cal X}_t)(m)dt=-z\int^1_0\alpha({\cal X}_t)(m)dt.$$
Hence the Hamiltonian vector field $X_{\widetilde F}$ of $\widetilde F$
with respect to $\widetilde\omega$ is given by  $\int^1_0\widetilde X_t dt$, and
the restriction of it to $(\partial M\times [1, +\infty), d\Theta)$ is given by
$$(m,z)\mapsto ({\cal X}(m,z), 0):=(\int^1_0{\cal X}_t(m)dt, 0).$$
This shows that the whole flow of $X_{\widetilde F}$ on $\widetilde M$, denoted by
$\phi^s_{\widetilde F}$, exists and on $\partial M\times [1,+\infty)$ has the form:
$\phi^s_{\widetilde F}(m,z)=(\chi^s(m), z)$, where $\chi^s$ is the flow of ${\cal X}$ on
$\partial M$. The key point is
$$\phi^s_{\widetilde F}(M)=M\quad{\rm and}\quad\phi^s_{\widetilde F}(\partial M\times (1, +\infty))=
\partial M\times (1, +\infty)\leqno(A.4)$$
for all $s\in\mbox{\Bb R}$.
Taking a strictly increasing smooth function $\eta:[0, 1/4]\to [0, 1]$ such that
$\eta(0)=0$, $\eta(1/4)=1$ and $\eta^\prime(1/4)=0$,  denoted by
$$\widetilde\psi_t:=\cases{ \widetilde\phi_{\eta(t)} & if $0\le t\le 1/4$,\cr
\widetilde\phi_1 &if $1/4\le t\le 3/4$,\cr \phi_{\widetilde
F}^{\eta(1-t)-1}\circ\widetilde\phi_1 & if $3/4\le t\le
1$.\cr}\leqno(A.5)$$ Setting $\widetilde
Z_t:=\frac{d}{dt}\widetilde\psi_t\circ\widetilde\psi_t^{-1}$, it
is a smooth family of vector fields on $\widetilde M$ and
$$\int^1_0\widetilde Z_t dt=0.\leqno(A.6)$$
>From (118) (A.4) (A.5) it follows that
$$\widetilde\psi_t(M)=M\quad{\rm and}\quad\widetilde\psi_t((\partial M\times (1, +\infty))=
\partial M\times (1, +\infty)\leqno(A.7)$$
for all $t\in [0,1]$. The straightforward computation shows that
$$\widetilde Z_t=\cases{\eta^\prime(t)\widetilde X_{\eta(t)} & if $0\le t\le 1/4$,\cr
0 &if $1/4\le t\le 3/4$,\cr
-\eta^\prime(1-t)X_{\widetilde F} & if $3/4\le t\le 1$;\cr}$$
$$\widetilde Z_t|_{\partial M\times [1,\infty)}(m,z)=\cases{\eta^\prime(t)({\cal X}_{\eta(t)}(m), 0) &
if $0\le t\le 1/4$,\cr 0 &if $1/4\le t\le 3/4$,\cr
-\eta^\prime(1-t)(\int^1_0{\cal X}_t(m)dt, 0) & if $3/4\le t\le
1$.\cr}\leqno(A.8)$$ Setting $\widetilde Y_t:=-\int^t_0\widetilde
Z_\lambda d\lambda$, then
$$\widetilde Y_t|_{\partial M\times [1,\infty)}(m,z)=\cases{(\int_0^{\eta(t)}{\cal X}_s(m)ds, 0) &
if $0\le t\le 1/4$,\cr (\int_0^1{\cal X}_s(m)ds, 0) &if $1/4\le
t\le 3/4$,\cr (\eta(1-t)\int^1_0{\cal X}_s(m)ds, 0) & if $3/4\le
t\le 1$.\cr}\leqno(A.9)$$ Let $\mbox{\Bb R}\to{\rm
Symp}_0(\widetilde
M,\widetilde\omega),\,s\mapsto\widetilde\theta^s_t$ be the flow
generated by $\widetilde Y_t$. Its existence is clear and is
uniquely determined by
$$\frac{d}{ds}\widetilde\theta^s_t=\widetilde Y_t\circ\widetilde\theta^s_t,\quad\widetilde\theta^0_t=id.$$
Moreover, since $\widetilde Y_0=\widetilde Y_1=0$ we get
$$\widetilde\theta^s_0=\widetilde\theta^s_1=id,\;\forall s\in\mbox{\Bb R}.$$
The key point is that $\widetilde\theta^s_t|_{\partial M\times [1,+\infty)}$ has the form
$$\widetilde\theta^s_t(m,z)=(\hat\theta^s_t(m), z)$$
for all $t\in [0, 1]$, $s\in\mbox{\Bb R}$ and $(m,z)\in\partial M\times [1, +\infty)$.
Here $\hat\theta^s_t:\partial M\to\partial M$.
 Setting $\widetilde\varphi_t:=\widetilde\theta^1_t\circ\widetilde\psi_t$
then it is easy to verify that ${\rm
Flux}(\{\widetilde\varphi_t\}_{0\le t\le T})=0$ for every $T\in
[0, 1]$. Thus it is an Hamiltonian path starting from $id$. Define
another Hamiltonian path starting from $id$, $[0, 1]\to {\rm
Ham}_0(\widetilde M,\widetilde\omega), \;t\mapsto\gamma_t$ by
$\gamma_t=id$ for $0\le t\le 3/4$, and $\gamma_t=\phi_{\widetilde
F}^{1-\eta(1-t)}$ for $3/4\le t\le 1$. Then
$t\mapsto\gamma_t\circ\widetilde\varphi_t$ is still an Hamiltonian
path starting from $id$. Moreover, when $s$ varies from $0$ to $1$
the path
$(\gamma_t\circ\widetilde\theta^s_t\circ\widetilde\psi_t)_{t\in
[0,1]}$ starting from $id$ varies from
$(\gamma_t\circ\widetilde\psi_t)_{t\in [0,1]}$ to
$(\gamma\circ\widetilde\varphi_t)_{t\in [0,1]}$ with fixed
endpoints. Since $\gamma_t(m,z)=(m,z)$ for $0\le t\le 3/4$, and
$\gamma_t(m,z)=(\chi^{1-\eta(1-t)}(m), z)$ for $3/4\le t\le 1$, it
is easily checked that
$$\gamma_t\circ\widetilde\theta^s_t\circ\widetilde\psi_t|_{\partial M\times [1, +\infty)}(m,z)
=\cases{(\hat\theta^s_t(\Phi_{\eta(t)}(m), z) & if $0\le t\le
1/4$,\cr (\hat\theta^s_t(\Phi_1(m)), z) &if $1/4\le t\le 3/4$,\cr
(\chi^{1-\eta(1-t)}\circ\hat\theta^s_t\circ{\chi}^{\eta(1-t)-1}\circ\Phi_1(m),
z) & if $3/4\le t\le 1$.\cr}\leqno(A.10)$$
>From these it follows that when $s$ varies from $0$ to $1$ the path
$(\gamma_t\circ\widetilde\theta^s_t\circ\widetilde\psi_t|_M)_{t\in [0,1]}$ varies from
$(\gamma_t\circ\widetilde\psi_t|_M)_{t\in [0,1]}$ to the Hamiltonian path
 $(\gamma_t\circ\widetilde\varphi_t|_M)_{t\in [0,1]}$
with fixed endpoints. But $\gamma_t\circ\widetilde\psi_t|_M=\phi_{\eta(t)}$ for
$0\le t\le 1/4$, and $\gamma_t\circ\widetilde\psi_t|_M=\phi_1$ for $1/4\le t\le 1$.
That is, $t\mapsto\gamma_t\circ\widetilde\psi_t|_M$ is only an
reparametrization of the path $t\mapsto\phi_t$.
This completes the proof of Theorem A. \hfill$\Box$

\end{document}